\documentclass[11pt]{article}
\usepackage{amsfonts}
\usepackage{amsfonts}
\usepackage{amsfonts}
\usepackage{mathrsfs}
\usepackage{bm}
\usepackage{cite}
\usepackage[colorlinks,linkcolor=blue,citecolor=blue]{hyperref}
\usepackage{amssymb,amsmath} 
 \allowdisplaybreaks

\catcode`!=11
\let\!int\int \def\int{\displaystyle\!int}
\let\!lim\lim \def\lim{\displaystyle\!lim}
\let\!sum\sum \def\sum{\displaystyle\!sum}
\let\!sup\sup \def\sup{\displaystyle\!sup}
\let\!inf\inf \def\inf{\displaystyle\!inf}
\let\!cap\cap \def\cap{\displaystyle\!cap}
\let\!max\max \def\max{\displaystyle\!max}
\let\!min\min \def\min{\displaystyle\!min}
\let\!frac\frac \def\frac{\displaystyle\!frac}
\catcode`!=12

\let\oldsection\section
\renewcommand\section{\setcounter{equation}{0}\oldsection}

\allowdisplaybreaks
\def\pf{\it{Proof.}\rm\quad}

\def\N{\mathbb{N}}

\newtheorem{thm}{Theorem}[section]
\newtheorem{lem}[thm]{Lemma}
\newtheorem{cor}[thm]{Corollary}

\begin{document}
\title {\bf Evaluations of Euler type sums of weight $\leq 5$}
\author{
{Ce Xu\thanks{Corresponding author. Email: 15959259051@163.com}}\\[1mm]
\small School of Mathematical Sciences, Xiamen University\\
\small Xiamen
361005, P.R. China}

\date{}
\maketitle\date{}
\maketitle \noindent{\bf Abstract } Let $p,p_1,\ldots,p_m$ be positive integers with $p_1\leq p_2\leq\cdots\leq p_m$ and $x\in [-1,1)$, define the so-called Euler type sums ${S_{{p_1}{p_2} \cdots {p_m},p}}\left( x \right)$, which are the infinite sums whose general term is a product of harmonic numbers of index $n$, a power of $n^{-1}$ and variable $x^n$, by
\[{S_{{p_1}{p_2} \cdots {p_m},p}}\left( x \right): = \sum\limits_{n = 1}^\infty  {\frac{{H_n^{\left( {{p_1}} \right)}H_n^{\left( {{p_2}} \right)} \cdots H_n^{\left( {{p_m}} \right)}}}{{{n^p}}}{x^n}}\quad (m\in \N:=\{1,2,3,\ldots\}),\]
where $H_n^{(p)}$ is defined by the generalized harmonic number.
Extending earlier work about classical Euler sums, we prove that whenever $p+p_1+\cdots+p_m \leq 5$, then all sums ${S_{{p_1}{p_2} \cdots {p_m},p}}\left( 1/2\right)$ can be expressed as a rational linear combination of products of zeta values, polylogarithms and $\log(2)$. The proof involves finding and solving linear equations which relate the different types of sums to each other.

\\[2mm]
\tableofcontents
\section{Introduction}
This paper is concerned with the discussion of sums of the type
\begin{align}\label{1.1}
{S_{{p_1}{p_2} \cdots {p_m},p}}\left( x \right): = \sum\limits_{n = 1}^\infty  {\frac{{H_n^{\left( {{p_1}} \right)}H_n^{\left( {{p_2}} \right)} \cdots H_n^{\left( {{p_m}} \right)}}}{{{n^p}}}{x^n}},
\end{align}
where the notation $H_n^{(p)}$ denotes the generalized harmonic number defined by (\cite{A1985,B1985})
\begin{align*}
H_n^{(p)}:=\sum\limits_{j=1}^n\frac {1}{j^p} \quad {\rm and} \quad H_n=H_n^{(1)}.
\end{align*}
For which values of the integer parameters $p,p_j\ (j=1,2,\ldots,m)$ and $x=\frac 1{2}$ can these sums be expressed in terms of the simpler values of polylogarithm function ${\rm Li}_p(x)$ and Riemann zeta function $\zeta(s)$? The polylogarithm function and Riemann zeta function are defined by (\cite{A2000})
\begin{align*}
&{\rm Li}{_p}\left( x \right) := \sum\limits_{n = 1}^\infty  {\frac{{{x^n}}}{{{n^p}}}}, \Re (p)>1,\\
&\zeta(s):=\sum\limits_{n = 1}^\infty {\frac {1}{n^{s}}},\Re(s)>1,
\end{align*}
with ${\rm Li_1}=-\log(1-x),\ x\in [-1,1).$
Here the quantity $w:={p _1} +  \cdots  + {p _m} + p$ is called the weight and the quantity $m$ is called the depth of ${S_{{p_1}{p_2} \cdots {p_m},p}}\left( x \right)$. As usual, repeated summands in partitions are indicated by powers, so that for instance
\[{S_{{1^2}{2^3}4,p}}\left(x\right) = {S_{112224,p}}\left(x\right) = \sum\limits_{n = 1}^\infty  {\frac{{H_n^2[H^{(2)} _n]^3{H^{(4)} _n}}}{{{n^p}}}}x^n. \]
When $x\rightarrow 1$, then the sum ${S_{{p_1}{p_2} \cdots {p_m},p}}\left( x \right)$ reduces to the classical Euler sum, which is defined by (\cite{FS1998,X2017,Xu2017})
\begin{align}\label{1.2}
{S_{{p_1}{p_2} \cdots {p_m},p}} := \sum\limits_{n = 1}^\infty  {\frac{{H_n^{\left( {{p_1}} \right)}H_n^{\left( {{p_2}} \right)} \cdots H_n^{\left( {{p_m}} \right)}}}
{{{n^p}}}}\quad (p>1),
\end{align}
which is also called the generalized (nonlinear) Euler sums.

Let $s_1,\ldots,s_k$ be positive integers. The multiple harmonic sums (MHS) are defined by (\cite{Xu2017})
\begin{align}\label{1.3}
{\zeta _n}\left( {{s_1},{s_2}, \cdots ,{s_k}} \right): = \sum\limits_{n \ge {n_1} > {n_2} >  \cdots  > {n_k} \ge 1} {\frac{1}{{n_1^{{s_1}}n_2^{{s_2}} \cdots n_k^{{s_k}}}}} ,
\end{align}
when $n<k$, then ${\zeta _n}\left( {{s_1},{s_2}, \cdots ,{s_k}} \right)=0$, and ${\zeta _n}\left(\emptyset \right)=1$. The integers $k$ and $w:=s_1+\ldots+s_k$ are called the depth and the weight of a multiple harmonic sum. For convenience, by ${\left\{ {{s_1}, \ldots ,{s_j}} \right\}_d}$ we denote the sequence of depth $dj$ with $d$ repetitions of ${\left\{ {{s_1}, \ldots ,{s_j}} \right\}}$. For example,
\[\zeta_n\left( {5,3,{{\left\{ 1 \right\}}_2}} \right) = \zeta_n\left( {5,3,1,1} \right),\;{\zeta_n }\left( {4,2,{{\left\{ 1 \right\}}_3}} \right) = {\zeta_n ^ \star }\left( {4,2,1,1,1} \right).\]
 When taking the limit $n\rightarrow \infty$ we get the so-called the multiple zeta value (MZV for short) (\cite{BBV2010,BBBL1997,Bro2013,DZ2012,Z2015}):
\begin{align}
&\zeta \left( {{s_1},{s_2}, \cdots ,{s_k}} \right) = \mathop {\lim }\limits_{n \to \infty } \zeta \left( {{s_1},{s_2}, \cdots ,{s_k}} \right),
\end{align}
defined for $s_2,\ldots,s_k\geq 1$ and $s_1\geq 2$ to ensure convergence of the series. It is obvious that $H_n^{(p)}=\zeta_n(p)$.

Moreover, we put a bar on top of $s_j\ (j=1,2,\cdots, k)$ if there is a sign $(-1)^{k_j}$ appearing in the denominator on the right of (\ref{1.3}), the sums are also called the alternating MHS. For example
\[{\zeta _n}\left( {{{\bar s}_1},{s_2}, \ldots ,{{\bar s}_k}} \right) = \sum\limits_{1 \le {k_m} <  \cdots  < {k_1} \le n} {\frac{{{{\left( { - 1} \right)}^{{k_1} + {k_m}}}}}{{k_1^{{s_1}} \cdots k_m^{{s_m}}}}},\]
\[{\zeta _n}\left( {\bar 2,3,\bar 1,4} \right)= \sum\limits_{n \ge {n_1} > {n_2} > {n_3} > {n_4} \ge 1}^{} {\frac{{{{\left( { - 1} \right)}^{{n_1} + {n_3}}}}}{{n_1^2n_2^3{n_3}n_4^4}}} .\]
The limit cases of alternating MHNs give rise to alternating multiple zeta values, for example
\begin{align*}
&{\zeta}\left( {\bar 2,3,5} \right) = \mathop {\lim }\limits_{n \to \infty } \zeta _n\left( {\bar 2,3,5} \right) = \sum\limits_{{n_1}> {n_2} >1}^{} {\frac{{{{\left( { - 1} \right)}^{{n_1}}}}}{{n_1^2n_2^3n_3^5}}} \\
&\zeta \left( {\bar 2,3,\bar 1,4} \right) = \mathop {\lim }\limits_{n \to \infty } {\zeta _n}\left( {\bar 2,3,\bar 1,4} \right) = \sum\limits_{{n_1} > {n_2} > {n_3} > {n_4} \ge 1}^{} {\frac{{{{\left( { - 1} \right)}^{{n_1} + {n_3}}}}}{{n_1^2n_2^3{n_3}n_4^4}}} .
\end{align*}

The relations between Euler sums and multiple zeta values have attracted a lot of research in the area in the last two decades. For details and historical introductions, please see \cite{BBG1994,BBV2010,BBG1995,BBBL1997,BBBL1999,Bro2013,FS1998,H1992,M2014,S2012,W2017,Xu2016,X2016,X2017,Xu2017,DZ2012,Z2015} and references therein. The origin of these numbers goes back to the correspondence of Euler with Goldbach in 1742-1743 (see \cite{BBG1994,FS1998}) that appeared in 1776. Euler's original contribution was a method to reduce double zeta values $\zeta(q,p)$ (or linear sums $S_{p,q}$) to certain rational linear combinations of products of zeta values, and established some important relation formula for them. For example, Euler proved that the linear sums $S_{p,q}$ are reducible to zeta values whenever $p+q$ is less than 7 or when $p+q$ is odd and less than 13, and he proved that
\[S_{1,k}=\sum\limits_{n = 1}^\infty  {\frac{{{H_n}}}{{{n^k}}}}  = \frac{1}{2}\left\{ {\left( {k + 2} \right)\zeta \left( {k + 1} \right) - \sum\limits_{i = 1}^{k - 2} {\zeta \left( {k - i} \right)\zeta \left( {i + 1} \right)} } \right\},\ k\geq2,\]
which, in particular implies the simplest but nontrivial relation
\begin{align*}
S_{1,2}=2\zeta(3)\quad {\rm or\ equivalently,} \quad \zeta(2,1)=\zeta(3).
\end{align*}
Moreover, he conjectured that the linear sums $S_{p,q}$ would be reducible whenever weight is odd, and even gave what he hoped to be the general formula. The conjecture was first proved by Borwein et al. \cite{BBG1995}. So, the linear sums $S_{p,q}$ can be evaluated in terms of zeta values in the following cases: $p=1,p=q,p+q$ odd and $p+q=6$ with $q\geq 2$\ (for more details, see \cite{BBG1994,BBG1995,FS1998}). Some examples on linear Euler sums follows
\begin{align*}
&{S_{3,2}} =  \frac{11}{2}\zeta \left( 5 \right) -2\zeta \left( 2 \right)\zeta \left( 3 \right),\\
&{S_{2,4}} = -\frac{{1}}{{3}}\zeta \left( 6 \right)+ {\zeta ^2}\left( 3 \right),\\
&{S_{5,2}} =  11\zeta \left( 7 \right)-4\zeta \left( 2 \right)\zeta \left( 5 \right) - 2\zeta \left( 3 \right)\zeta \left( 4 \right) ,\\
&{S_{4,3}} = -17\zeta \left( 7 \right) +\zeta \left( 3 \right)\zeta \left( 4 \right)+ 10\zeta \left( 2 \right)\zeta \left( 5 \right),\\
 &{S_{4,5}} =  - \frac{{125}}{2}\zeta \left( 9 \right) + 35\zeta \left( 2 \right)\zeta \left( 7 \right) + 5\zeta \left( 4 \right)\zeta \left( 5 \right), \\
 &{S_{2,7}} =  - \frac{{35}}{2}\zeta \left( 9 \right) + 7\zeta \left( 2 \right)\zeta \left( 7 \right) + 2\zeta \left( 3 \right)\zeta \left( 6 \right) + 4\zeta \left( 4 \right)\zeta \left( 5 \right), \\
& {S_{4,7}} =  - \frac{{229}}{2}\zeta \left( {11} \right) + 84\zeta \left( 2 \right)\zeta \left( 9 \right) + 21\zeta \left( 4 \right)\zeta \left( 7 \right) + 4\zeta \left( 5 \right)\zeta \left( 6 \right).
\end{align*}

Investigation of Euler sums has a long history, but usually the authors were not aware of Euler's results, so that special instances of Euler's identities have been independently rediscovered time and again. It was mainly the publication of B. Berndt's edition of Ramanujan's notebooks \cite{B1985} that served to fit all the scattered individual results into the framework of Euler's work. Besides the works referred to above, there are many other researches devoted to the Euler sums. For example, in 1994, Bailey et al. \cite{BBG1994} proved that all Euler sums of the form $S_{1^p,q}$ for weights $p+q\in \{3,4,5,6,7,9\}$ are reducible to {\bf Q}-linear combinations of zeta values by using the experimental method. In 1995, Borwein et al. \cite{BBG1995} showed that the quadratic sums $S_{ 1^2,q}$ can reduce to linear sums $S_{2,q}$ and polynomials in zeta values. In 1998, Flajolet and Salvy \cite{FS1998} used the contour integral representations and residue computation to show that the quadratic sums $S_{p_1p_2,q}$ are reducible to linear sums and zeta values when the weight $p_1 + p_2 + q$ is even and $p_1,p_2>1$. The best results to date are due to Xu and Wang et al, see the most recent papers \cite{W2017,X2017,Xu2017}. In \cite{X2017,Xu2017}, we proved that all Euler sums of weight $\leq 8$ are reducible to $\mathbb{Q}$-linear combinations of single zeta monomials with the addition of $\{S_{2,6}\}$ for weight 8. For weight 9, all Euler sums of the form ${S_{{s_1} \cdots {s_k},q}}$ with $q\in \{4,5,6,7\}$ are expressible polynomially in terms of zeta values. For weight $p_1+p_2+q=10$, all quadratic sums $S_{p_1p_2,q}$ are reducible to $S_{2,6}$ and $S_{2,8}$.
Wang et al \cite{W2017} shown that all Euler sums of weight $\leq 9$ are reducible to zeta values and linear sums. Examples for such evaluations, all due to Xu and Wang, are
\begin{align*}
&S_{12^2,3} =  - \frac{{6313}}{{288}}\zeta \left( 8 \right) + \frac{{43}}{2}\zeta \left( 3 \right)\zeta \left( 5 \right) + \frac{1}{2}\zeta \left( 2 \right){\zeta ^2}\left( 3 \right) - \frac{{17}}{4}{S_{2,6}},\\
&{S_{{1^2}23,2}} = \frac{{505}}{{36}}\zeta (9) + \frac{7}{4}\zeta (2)\zeta (7) + 3\zeta (3)\zeta (6) - \frac{{37}}{4}\zeta (4)\zeta (5) - \frac{5}{3}{\zeta ^3}(3),\\
&{S_{{{12}^2},4}}  =  - \frac{{775}}{{36}}\zeta \left( 9 \right) + \frac{{85}}{8}\zeta \left( 2 \right)\zeta \left( 7 \right) - \frac{{221}}{{24}}\zeta \left( 3 \right)\zeta \left( 6 \right) + 10\zeta \left( 4 \right)\zeta \left( 5 \right) + 3{\zeta ^3}\left( 3 \right).
\end{align*}

In this paper we are interested in Euler-type sums with hyperharmonic numbers ${S_{{p_1}{p_2} \cdots {p_m},p}}\left( x \right)$. Such series could be of interest in analytic number theory. We will show that these sums are related to the values of the Riemann zeta function and polylogarithm function when $x=\frac 1{2}$ and weight $\leq 5$ by using the method of based on simple integral representations of logarithms.
\section{Some lemmas}
In this section, we give some lemmas which will be useful in the development of the main results.
\begin{lem}(\cite{X2016})\label{lem 2.1}
Let $s,t$ be positive integers with $x\in [-1,1)$. Then the product of two polylogarithm functions are reducible to Euler type sums \begin{align}\label{2.1}
{\rm Li}{_s}\left( x \right){\rm Li}{_t}\left( x \right) =& \sum\limits_{j = 1}^s {A_j^{\left( {s,t} \right)}} \sum\limits_{n = 1}^\infty  {\frac{{{H^{(j)} _n}}}{{{n^{s + t - j}}}}} {x^n} + \sum\limits_{j = 1}^t {B_j^{\left( {s,t} \right)}} \sum\limits_{n = 1}^\infty  {\frac{{{H^{(j)} _n}}}{{{n^{s + t - j}}}}} {x^n}\nonumber \\& - \left( {\sum\limits_{j = 1}^s {A_j^{\left( {s,t} \right)}}  + \sum\limits_{j = 1}^t {B_j^{\left( {s,t} \right)}} } \right){\rm Li}{_{s + t}}\left( x \right),
\end{align}
where $A_j^{\left( {s,t} \right)} = \left( {\begin{array}{*{20}{c}}
   {s + t - j - 1}  \\
   {s - j}  \\
\end{array}} \right),B_j^{\left( {s,t} \right)} = \left( {\begin{array}{*{20}{c}}
   {s + t - j - 1}  \\
   {t - j}  \\
\end{array}} \right).$
\end{lem}

\begin{lem}(\cite{CX2016})\label{lem 2.2}
For integers $m\geq 1$ and $k\geq 0$, then the following identity holds:
\begin{align}\label{2.2}
&m!\sum\limits_{l = 0}^k {l!\left( {\begin{array}{*{20}{c}}
   k  \\
   l  \\
\end{array}} \right){{\left( {\log\left( 2\right)} \right)}^{k - l}}\zeta \left( {l + 2,{{\left\{ 1 \right\}}_{m - 1}};\frac{1}{2}} \right)}\nonumber \\
& + k!\sum\limits_{l = 0}^m {l!\left( {\begin{array}{*{20}{c}}
   m  \\
   l  \\
\end{array}} \right){{\left( {\log \left( 2\right)} \right)}^{m - l}}\zeta \left( {l + 1,{{\left\{ 1 \right\}}_{k }};\frac{1}{2}} \right)} \nonumber\\
&= m!k!\zeta \left( {m + 1,{{\left\{ 1 \right\}}_k}} \right),
\end{align}
where the $\zeta \left( {{s_1},{s_2}, \cdots ,{s_m};x} \right): = {\rm{L}}{{\rm{i}}_{{s_1},{s_2}, \cdots ,{s_m}}}\left( x \right)$ denotes the multiple polylogarithm function defined by
\begin{align}\label{2.3}
\zeta \left( {{s_1},{s_2}, \cdots ,{s_m};x} \right): ={\rm{L}}{{\rm{i}}_{{s_1},{s_2}, \cdots ,{s_m}}}\left( x \right): = \sum\limits_{1 \le {k_m} <  \cdots  < {k_1}} {\frac{{{x^{{k_1}}}}}{{k_1^{{s_1}}k_2^{{s_2}} \cdots k_m^{{s_m}}}}} ,\;x \in \left[ { - 1,1} \right).
\end{align}
Here ${\bf S}:=(s_1,s_2,\ldots,s_m)\in (\N)^m$ in above definition (\ref{2.3}). Of course, if $s_1>1$, then we can allow $x=1$.
\end{lem}

\begin{lem}(\cite{CX2017})\label{lem 2.3}
For integer $m\in\N_0:=\N\cup \{0\}$, then the following identity holds:
\begin{align}\label{2.4}
\zeta\left(2,\{1\}_m;\frac 1{2}\right)={\rm{L}}{{\rm{i}}_{2,{{\left\{ 1 \right\}}_m}}}\left( {\frac{1}{2}} \right) = \zeta \left( {m + 2} \right) - \sum\limits_{l = 0}^{m + 1} {\frac{{{{\left( {\log( 2)} \right)}^{m + 1 - l}}}}{{\left( {m + 1 - l} \right)!}}{\rm{L}}{{\rm{i}}_{l + 1}}\left( {\frac{1}{2}} \right)} .
\end{align}
\end{lem}

\begin{lem}(\cite{Xu2017})\label{lem 2.4}
For integer $k>0$ and $x\in [-1,1)$, then we have
\begin{align}\label{2.5}
{\log ^k}\left( {1 - x} \right) = {\left( { - 1} \right)^k}k!\sum\limits_{n = 1}^\infty  {\frac{{{x^n}}}{n}{\zeta _{n - 1}}\left( {{{\left\{ 1 \right\}}_{k - 1}}} \right)},\end{align}

\begin{align}\label{2.6}
s\left( {n,k} \right) = \left( {n - 1} \right)!{\zeta _{n - 1}}\left( {{{\left\{ 1 \right\}}_{k - 1}}} \right),
\end{align}
where ${s\left( {n,k} \right)}$ denotes the (unsigned) Stirling number of the first kind (see \cite{L1974}), and we have
\begin{align*}
& s\left( {n,1} \right) = \left( {n - 1} \right)!,\\&s\left( {n,2} \right) = \left( {n - 1} \right)!{H_{n - 1}},\\&s\left( {n,3} \right) = \frac{{\left( {n - 1} \right)!}}{2}\left[ {H_{n - 1}^2 - {H^{(2)} _{n - 1}}} \right],\\
&s\left( {n,4} \right) = \frac{{\left( {n - 1} \right)!}}{6}\left[ {H_{n - 1}^3 - 3{H_{n - 1}}{H^{(2)} _{n - 1}} + 2{H^{(3)}_{n - 1}}} \right], \\
&s\left( {n,5} \right) = \frac{{\left( {n - 1} \right)!}}{{24}}\left[ {H_{n - 1}^4 - 6{H^{(4)}_{n - 1}} - 6H_{n - 1}^2{H^{(2)}_{n - 1}} + 3(H^{(2)}_{n-1})^2+ 8H_{n - 1}^{}{H^{(3)}_{n - 1}}} \right].
\end{align*}
The Stirling numbers ${s\left( {n,k} \right)}$ of the first kind satisfy a recurrence relation in the form
\[s\left( {n,k} \right) = s\left( {n - 1,k - 1} \right) + \left( {n - 1} \right)s\left( {n - 1,k} \right),\;\;n,k \in \N,\]
with $s\left( {n,k} \right) = 0,n < k,s\left( {n,0} \right) = s\left( {0,k} \right) = 0,s\left( {0,0} \right) = 1$.
\end{lem}

\begin{lem}(\cite{X2016})\label{lem 2.5}
For integers $n\geq1$ and $ k\geq 0$, then
\begin{align}\label{2.7}
\int\limits_0^1 {{t^{n - 1}}{{\log }^k}\left( {1 - t} \right)} dt = {\left( { - 1} \right)^k}\frac{{{Y_k}\left( n \right)}}{n},\ k,n\in \N,
\end{align}
where ${Y_k}\left( n \right) = {Y_k}\left( {{H _n}( 1 ),1!{H^{(2)} _n},2!{H^{(3)} _n}, \cdots ,
\left( {r - 1} \right)!{H^{(r)} _n}, \cdots } \right)$, ${Y_k}\left( {{x_1},{x_2}, \cdots } \right)$ stands for the complete exponential Bell polynomial is defined by (see \cite{L1974})
\[\exp \left( {\sum\limits_{m \ge 1}^{} {{x_m}\frac{{{t^m}}}{{m!}}} } \right) = 1 + \sum\limits_{k \ge 1}^{} {{Y_k}\left( {{x_1},{x_2}, \cdots } \right)\frac{{{t^k}}}{{k!}}}.\]
 From the definition of the complete exponential Bell polynomial, we deduce
$${Y_1}\left( n \right) = {H_n},{Y_2}\left( n \right) = H_n^2 + {H^{(2)} _n},{Y_3}\left( n \right) =  H_n^3+ 3{H_n}{H^{(2)} _n}+ 2{H^{(3)} _n},$$
\[{Y_4}\left( n \right) = H_n^4 + 8{H_n}{H^{(3)} _n} + 6H_n^2{H^{(2)} _n} + 3(H^{(2)} _n)^2 + 6{H^{(4)} _n},\]
\[{Y_5}\left( n \right) = H_n^5 + 10H_n^3{H^{(2)} _n} + 20H_n^2{H^{(3)}_n} + 15{H_n}({H^{(2)}_n})^2 + 30{H_n}{H^{(4)} _n}+ 20{H^{(2)} _n}{H^{(3)} _n} + 24{H^{(5)} _n},\]
Moreover, from \cite{Xu2017}, we know that the ${Y_k}\left( n \right)$ is a rational linear combination of products of harmonic numbers.
\end{lem}

\section{Main theorems and proofs}
In this section, we will establish some relations between Euler type sums ${S_{{p_1}{p_2} \cdots {p_m},p}}\left( 1/2\right)$ and integrals of logarithms by using above lemmas.
\begin{thm} \label{thm3.1}
For any $m\in \N_0$, then the following identity holds:
\begin{align}\label{3.1}
S_{1,m+1}(-1)=&  \zeta \left( \overline{m + 2} \right) - \frac{{{{\left( { - 1} \right)}^{m + 1}}}}{{m!\left( {m + 2} \right)}}{\left( {\log (2)} \right)^{m + 2}}\nonumber\\
&- \frac{{{{1}}}}{{m!}}\sum\limits_{k = 1}^m {\left( {\begin{array}{*{20}{c}}
   m  \\
   k  \\
\end{array}} \right)} {\left( { - 1} \right)^{ k}}\int\limits_{1/2}^1 {\frac{{{{\log }^{m - k + 1}}(x){{\log }^k}\left( {1 - x} \right)}}{x}dx}.
\end{align}
\end{thm}
\pf By a direct calculation, we deduce that
\begin{align}\label{3.2}
{S_{1,m + 1}}\left( { - 1} \right)& = \sum\limits_{n = 1}^\infty  {\frac{{{H_n}}}{{{n^{m + 1}}}}{{\left( { - 1} \right)}^n}}  = \frac{{{{\left( { - 1} \right)}^m}}}{{m!}}\sum\limits_{n = 1}^\infty  {{{\left( { - 1} \right)}^n}{H_n}\int\limits_0^1 {{x^{n - 1}}{{\log }^m}\left( x \right)dx} }\nonumber \\
&= \frac{{{{\left( { - 1} \right)}^{m + 1}}}}{{m!}}\int\limits_0^1 {\frac{{{{\log }^m}\left( x \right)\log \left( {1 + x} \right)}}{{\left( {1 + x} \right)x}}} dx\nonumber\\
& = \frac{{{{\left( { - 1} \right)}^{m + 1}}}}{{m!}}\left\{ {\int\limits_0^1 {\frac{{{{\log }^m}\left( x \right)\log \left( {1 + x} \right)}}{x}} dx - \int\limits_0^1 {\frac{{{{\log }^m}\left( x \right)\log \left( {1 + x} \right)}}{{1 + x}}} dx} \right\}\nonumber\\
&=\frac{{{{\left( { - 1} \right)}^m}}}{{m!}}\sum\limits_{n = 1}^\infty  {\frac{{{{\left( { - 1} \right)}^n}}}{n}\int\limits_0^1 {{x^{n - 1}}{{\log }^m}\left( x \right)dx} }  + \frac{{{{\left( { - 1} \right)}^m}}}{{m!}}\int\limits_1^2 {\frac{{{{\log }^m}\left( {t - 1} \right)\log \left( t \right)}}{t}dt} \nonumber\\
& = \zeta \left( \overline {m + 2} \right) - \frac{{{{\left( { - 1} \right)}^m}}}{{m!}}\int\limits_{1/2}^1 {\frac{{{{\log }^m}\left( {\frac{{1 - x}}{x}} \right)\log \left( x \right)}}{x}dx}.
\end{align}
We note that the integral on the right hand side of (\ref{3.2}) can be rewritten as
\begin{align}\label{3.3}
&\int\limits_{1/2}^1 {\frac{{{{\log }^m}\left( {\frac{{1 - x}}{x}} \right)\log \left( x \right)}}{x}dx}  = \sum\limits_{k = 0}^m {\left( {\begin{array}{*{20}{c}}
   m  \\
   k  \\
\end{array}} \right){{\left( { - 1} \right)}^{m - k}}\int\limits_{1/2}^1 {\frac{{{{\log }^k}\left( {1 - x} \right){{\log }^{m - k + 1}}\left( x \right)}}{x}dx} } \nonumber\\
 =& \sum\limits_{k = 1}^m {\left( {\begin{array}{*{20}{c}}
   m  \\
   k  \\
\end{array}} \right){{\left( { - 1} \right)}^{m - k}}\int\limits_{1/2}^1 {\frac{{{{\log }^k}\left( {1 - x} \right){{\log }^{m - k + 1}}\left( x \right)}}{x}dx} }  - \frac{1}{{m + 2}}{\log ^{m + 2}}\left( 2 \right).
\end{align}
Then, substituting identity (\ref{3.3}) into (\ref{3.2}) yields the desired result. This completes the proof of Theorem \ref{thm3.1}. \hfill$\square$

\begin{thm}\label{thm3.2} For integer $m\in \N_0$ and real $x\in [-1,1)$, we have
\begin{align}\label{3.4}
\int\limits_0^x {\frac{{{{\log }^m}\left( {1 - t} \right)}}{{1 + t}}dt}  =& {\left( { - 1} \right)^m}m!{\rm{L}}{{\rm{i}}_{m + 1}}\left( {\frac{1}{2}} \right) + {\log ^m}\left( {1 - x} \right)\log \left( {\frac{{1 + x}}{2}} \right)\nonumber
\\
&+ \sum\limits_{l = 1}^m {{{\left( { - 1} \right)}^{l + 1}}{{\log }^{m - l}}\left( {1 - x} \right){{\left( m \right)}_l}{\rm Li}{_{l + 1}}\left( {\frac{{1 - x}}{2}} \right)},
\end{align}
where $(m)_l:=m(m-1)\cdots(m-l+1)$.
\end{thm}
\pf Changing the variable $t\mapsto 1-u$, then the integral on the left hand side of (\ref{3.4}) can be rewritten as
\begin{align}\label{3.5}
\int\limits_0^x {\frac{{{{\log }^m}\left( {1 - t} \right)}}{{1 + t}}dt}  = \int_{1 - x}^1 {\frac{{{{\log }^m}\left( u \right)}}{{2 - u}}} du = \sum\limits_{n = 1}^\infty  {\frac{1}{{{2^n}}}\int_{1 - x}^1 {{u^{n - 1}}{{\log }^m}\left( u \right)} } dt.\end{align}
On the other hand, by using integration by parts, we deduce that, for $n,m\in \N$,
\begin{align}\label{3.6}
\int\limits_0^x {{t^{n - 1}}{{\left( {\log (t)} \right)}^m}} dt = \sum\limits_{l = 0}^m {l!\left( {\begin{array}{*{20}{c}}
   m  \\
   l  \\
\end{array}} \right)\frac{{{{\left( { - 1} \right)}^l}}}{{{n^{l + 1}}}}{{\left( {\log (x)} \right)}^{m - l}}{x^n}},x\in (0,1).\end{align}
Hence, substituting (\ref{3.6}) into (\ref{3.5}), by a simple calculation we obtain the formula (\ref{3.4}).\hfill$\square$\\
Taking $m=1$ and $2$ in (\ref{3.4}), we obtain the following cases
\begin{align}\label{3.7}
\int\limits_0^x {\frac{{\log \left( {1 - t} \right)}}{{1 + t}}dt}  = \log \left( {1 - x} \right)\log \left( {\frac{{1 + x}}{2}} \right) + {\rm{L}}{{\rm{i}}_2}\left( {\frac{{1 - x}}{2}} \right) - {\rm{L}}{{\rm{i}}_2}\left( {\frac{1}{2}} \right),
\end{align}
\begin{align}\label{3.8}
\int\limits_0^x {\frac{{{{\log }^2}\left( {1 - t} \right)}}{{1 + t}}dt}  =& 2{\rm{L}}{{\rm{i}}_3}\left( {\frac{1}{2}} \right) + {\log ^2}\left( {1 - x} \right)\log \left( {\frac{{1 + x}}{2}} \right) \nonumber\\&+ 2\log \left( {1 - x} \right){\rm{L}}{{\rm{i}}_2}\left( {\frac{{1 - x}}{2}} \right) - 2{\rm{L}}{{\rm{i}}_3}\left( {\frac{{1 - x}}{2}} \right).
\end{align}
Noting that from \cite{BD2006}, we have
\[{\rm{L}}{{\rm{i}}_2}\left( {\frac{1}{2}} \right) = \frac{{\zeta \left( 2 \right) - {{\log }^2}\left( 2 \right)}}{2},{\rm{L}}{{\rm{i}}_3}\left( {\frac{1}{2}} \right) = \frac{7}{8}\zeta \left( 3 \right) - \frac{1}{2}\zeta \left( 2 \right)\log \left( 2 \right) + \frac{1}{6}{\log ^3}\left( 2 \right).\]
\begin{cor}\label{cor3.3} For any $x\in [-1,1)$, then the following identity holds:
\begin{align}\label{3.9}\sum\limits_{n = 1}^\infty  {\frac{{H_n^{\left( 2 \right)} - L_n^2\left( 1 \right)}}{{n + 1}}} {x^{n + 1}} = 4\left[ {{\rm{L}}{{\rm{i}}_3}\left( {\frac{{1 - x}}{2}} \right) - {\rm{L}}{{\rm{i}}_3}\left( {\frac{1}{2}} \right)} \right] - 2\log \left( {1 - x} \right)\left[ {{\rm{L}}{{\rm{i}}_2}\left( {\frac{{1 - x}}{2}} \right) + {\rm{L}}{{\rm{i}}_2}\left( {\frac{1}{2}} \right)} \right],\end{align}
where $L_n(p)$ denotes the alternating harmonic number, which is defined by
$$L_{n}(p):=\sum\limits_{j=1}^n\frac{(-1)^{j-1}}{j^p},\ p \in \N.$$
\end{cor}
\pf In \cite{XC2016}, we proved the result
\[\sum\limits_{n = 1}^\infty  {\frac{{H_n^{\left( 2 \right)} - L_n^2\left( 1 \right)}}{{n + 1}}} {x^{n + 1}} = 2\log \left( {1 - x} \right)\int\limits_0^x {\frac{{\log \left( {1 - t} \right)}}{{1 + t}}dt}  - 2\int\limits_0^x {\frac{{{{\log }^2}\left( {1 - t} \right)}}{{1 + t}}dt} .\]
Applying formulas (\ref{3.7}) and (\ref{3.8}) to the above equation, we deduce the desired result.\hfill$\square$
Differentiate both sides of (\ref{3.9}), then
\[\sum\limits_{n = 1}^\infty  {\left\{ {H_n^{\left( 2 \right)} - L_n^2\left( 1 \right)} \right\}{x^n}}  = 2\frac{{{\rm{L}}{{\rm{i}}_2}\left( {\frac{1}{2}} \right) - {\rm{L}}{{\rm{i}}_2}\left( {\frac{{1 - x}}{2}} \right) - \log \left( {1 - x} \right)\log \left( {\frac{{1 + x}}{2}} \right)}}{{1 - x}},\ x\in(-1,1).\]
\begin{thm}\label{thm3.4} For integers $m\in \N_0$ and $k\in \N$, then the following relations hold:
\begin{align}\label{3.10}
\sum\limits_{n = 1}^\infty  {\frac{{s\left( {n,m + 1} \right){Y_k}\left( n \right)}}{{n!{2^n}}}}  = \frac{{{{\left( { - 1} \right)}^k}}}{{m!}}\sum\limits_{j = 0}^m {{{\left( { - 1} \right)}^j}\left( {\begin{array}{*{20}{c}}
   m  \\
   j  \\
\end{array}} \right){{\log }^{m - j}}\left( 2 \right)\int\limits_0^1 {\frac{{{{\log }^j}\left( {1 + t} \right){{\log }^k}\left( t \right)}}{{1 + t}}dt} } ,
\end{align}
\begin{align}\label{3.11}
\int\limits_0^1 {\frac{{{{\log }^j}\left( {1 + t} \right){{\log }^k}\left( t \right)}}{{1 + t}}dt}  = {\left( { - 1} \right)^{k + j}}k!j!\sum\limits_{n = 1}^\infty  {\frac{{s\left( {n,j + 1} \right)}}{{n!{n^k}}}{{\left( { - 1} \right)}^{n + 1}}} .
\end{align}
\end{thm}
\pf The identity (\ref{3.11}) is easily derived. Next, we prove the formula (\ref{3.10}).
By using Lemma \ref{lem 2.4},\ \ref{lem 2.5} and considering the following integral
\[I(m,k):=\int\limits_0^1 {\frac{{{{\log }^m}\left( {1 - \frac{x}{2}} \right){{\log }^k}\left( {1 - x} \right)}}{{1 - \frac{x}{2}}}dx} ,\]
we have
\begin{align}\label{3.12}
I\left( {m,k} \right) =& {\left( { - 1} \right)^m}m!\sum\limits_{n = 1}^\infty  {\frac{{s\left( {n + 1,m + 1} \right)}}{{n!{2^n}}}\int\limits_0^1 {{x^n}{{\log }^k}\left( {1 - x} \right)dx} }\nonumber \\
= &{\left( { - 1} \right)^{m + k}}m!\sum\limits_{n = 1}^\infty  {\frac{{s\left( {n + 1,m + 1} \right){Y_k}\left( {n + 1} \right)}}{{n!{2^n}\left( {n + 1} \right)}}}\nonumber \\
= &2{\left( { - 1} \right)^{m + k}}m!\sum\limits_{n = 1}^\infty  {\frac{{s\left( {n,m + 1} \right){Y_k}\left( n \right)}}{{n!{2^n}}}} .
\end{align}
On the other hand, applying the change of variable $x\rightarrow 1-t$ to the above integral on the left hand side of (\ref{3.12}), which can be rewritten as
\begin{align}\label{3.13}
 I\left( {m,k} \right) =& 2\int\limits_0^1 {\frac{{{{\log }^m}\left( \frac {1 + t}{2} \right){{\log }^k}\left( t \right)}}{{1 + t}}dt}\nonumber  \\
  = &2\sum\limits_{j = 0}^m {{{\left( { - 1} \right)}^{m - j}}\left( {\begin{array}{*{20}{c}}
   m  \\
   j  \\
\end{array}} \right){{\log }^{m - j}}\left( 2 \right)\int\limits_0^1 {\frac{{{{\log }^j}\left( {1 + t} \right){{\log }^k}\left( t \right)}}{{1 + t}}dt} }
\end{align}
Thus, combining identities (\ref{3.12}) and (\ref{3.13}) we obtain the result. \hfill$\square$

\begin{thm}\label{thm3.5} For positive integers $m$ and $k$, then the following equation holds:
\begin{align}\label{3.14}
{\left( { - 1} \right)^{m + k}}m!\sum\limits_{n = 1}^\infty  {\frac{{s\left( {n,m} \right){Y_k}\left( n \right)}}{{n!n{2^n}}}}  = \sum\limits_{j = 0}^m {{{\left( { - 1} \right)}^j}\left( {\begin{array}{*{20}{c}}
   m  \\
   j  \\
\end{array}} \right){{\log }^j}\left( 2 \right)\int\limits_0^1 {\frac{{{{\log }^{m - j}}\left( {1 + t} \right){{\log }^k}\left( t \right)}}{{1 - t}}dt} }
\end{align}
\end{thm}
\pf Similarly as in the proof of Theorem \ref{thm3.4}, considering the integral
\[J(m,k):=\int\limits_0^1 {\frac{{{{\log }^m}\left( {1 - \frac{x}{2}} \right){{\log }^k}\left( {1 - x} \right)}}{x}dx}. \]
Then using the Lemma \ref{lem 2.4},\ \ref{lem 2.5} and applying the change of variable $x\rightarrow 1-t$ to the above integral, we have
\begin{align}\label{3.15}
 J\left( {m,k} \right) &= {\left( { - 1} \right)^m}m!\sum\limits_{n = 1}^\infty  {\frac{{s\left( {n,m} \right)}}{{n!{2^n}}}\int\limits_0^1 {{x^{n - 1}}{{\log }^k}\left( {1 - x} \right)dx} }\nonumber  \\
  &= {\left( { - 1} \right)^{m + k}}m!\sum\limits_{n = 1}^\infty  {\frac{{s\left( {n,m} \right){Y_k}\left( n \right)}}{{n!n{2^n}}}} \nonumber \\
 & = \int\limits_0^1 {\frac{{{{\log }^m}\left( {\frac{{1 + t}}{2}} \right){{\log }^k}\left( t \right)}}{{1 - t}}dt} \nonumber \\
 & = \sum\limits_{j = 0}^m {{{\left( { - 1} \right)}^j}\left( {\begin{array}{*{20}{c}}
   m  \\
   j  \\
\end{array}} \right){{\log }^j}\left( 2 \right)\int\limits_0^1 {\frac{{{{\log }^{m - j}}\left( {1 + t} \right){{\log }^k}\left( t \right)}}{{1 - t}}dt} }.
\end{align}
Thus, the formula (\ref{3.14}) holds. \hfill$\square$
\begin{thm}\label{thm3.6} For positive integer $m$, then the following identity holds:
\begin{align}\label{3.16}
\sum\limits_{n = 1}^\infty  {\frac{{{H_n}H_n^{\left( {m + 1} \right)}}}{{n{2^n}}}}  =& \sum\limits_{n = 1}^\infty  {\frac{{H_n^{\left( {m + 1} \right)}}}{{{n^2}{2^n}}}}  + \frac{{{{\left( { - 1} \right)}^{m + 1}}}}{{2\left( {m!} \right)}}\int\limits_0^1 {\frac{{{{\log }^2}\left( {1 + t} \right){{\log }^m}\left( {1 - t} \right)}}{t}dt}\nonumber \\
& - \frac{{{{\left( { - 1} \right)}^{m + 1}}}}{{m!}}\log \left( 2 \right)\int\limits_0^1 {\frac{{\log \left( {1 + t} \right){{\log }^m}\left( {1 - t} \right)}}{t}dt} .
\end{align}
\end{thm}
\pf In the same way as in proofs of Theorems \ref{thm3.4} and \ref{thm3.5}, considering the integral
\[\int\limits_0^1 {\frac{{{{\log }^2}\left( {1 - \frac{x}{2}} \right){{\log }^m}\left( x \right)}}{{1 - x}}dx}.\]
Then applying the identity (\ref{2.5}) with the help of the following elementary integral
\begin{align}\label{3.17}\int\limits_0^1 {\frac{{{x^n}{{\log }^m}x}}{{1 - x}}} dx = {\left( { - 1} \right)^m}m!\left( {\zeta \left( {m + 1} \right) - {H^{(m+1)} _n}} \right)\;\ (m \in \N),\end{align}
By a simple calculation, we deduce the result.\hfill$\square$
\section{Some results on integral of logarithms}
In \cite{Xu2016,Xyz2016,Xu2015,XC2016}, we obtain numerous results of some alternating Euler sums of weight $\leq 6$. We can use these results to find some nice evaluations of integral of logarithms. Hence, in this section, we will give many closed form representations of logarithms' integrals. By using Lemma \ref{lem 2.4},\ \ref{lem 2.5} and formula (\ref{3.6}) with the help of results of references \cite{Xu2016,Xyz2016,Xu2015,XC2016}, the following identities are easily derived
\begin{align}\label{4.1}
\int\limits_0^1 {\frac{{{{\log }^m}\left( {1 - x} \right)}}{x}dx = {{\left( { - 1} \right)}^m}m!\zeta \left( {m + 1} \right)} ,
\end{align}
\begin{align}\label{4.2}
\int\limits_0^1 {\frac{{{{\log }^m}\left( {1 - x} \right)}}{{1 - \frac{x}{2}}}} dx = 2{\left( { - 1} \right)^m}m!\left( {1 - \frac{1}{{{2^m}}}} \right)\zeta \left( {m + 1} \right),
\end{align}
\begin{align}\label{4.3}
 \int\limits_{1/2}^1 {\frac{{\log \left( x \right){{\log }^2}\left( {1 - x} \right)}}{x}dx}  =&  - 2{\rm{L}}{{\rm{i}}_4}\left( {\frac{1}{2}} \right) - \frac{1}{2}\zeta \left( 4 \right) + \frac{1}{4}\zeta \left( 3 \right)\log \left( 2 \right) \nonumber\\
  &- \frac{1}{3}{\log ^4}\left( 2 \right) + 2\sum\limits_{n = 1}^\infty  {\frac{{{H_n}}}{{{n^3}{2^n}}}} ,
\end{align}
\begin{align}\label{4.4}
 \int\limits_{1/2}^1 {\frac{{{{\log }^2}\left( x \right)\log \left( {1 - x} \right)}}{x}dx}  =& 2{\rm{L}}{{\rm{i}}_4}\left( {\frac{1}{2}} \right) - 2\zeta \left( 4 \right) + \frac{7}{4}\zeta \left( 3 \right)\log \left( 2 \right)\nonumber \\
  &- \frac{1}{2}\zeta \left( 2 \right){\log ^2}\left( 2 \right) - \frac{1}{6}{\log ^4}\left( 2 \right),
\end{align}
\begin{align}\label{4.5}
 \int\limits_0^1 {\frac{{\log \left( {1 + x} \right){{\log }^2}\left( {1 - x} \right)}}{x}dx}  = & 2{\rm{L}}{{\rm{i}}_4}\left( {\frac{1}{2}} \right) - \frac{5}{8}\zeta \left( 4 \right) + \frac{7}{4}\zeta \left( 3 \right)\log \left( 2 \right)\nonumber \\
 & - \frac{1}{2}\zeta \left( 2 \right){\log ^2}\left( 2 \right) + \frac{1}{{12}}{\log ^4}\left( 2 \right),
\end{align}
\begin{align}\label{4.6}
\int\limits_0^1 {\frac{{\log \left( {1 - \frac{x}{2}} \right){{\log }^3}\left( {1 - x} \right)}}{x}dx}  = 12\zeta \left( 5 \right) - \frac{{21}}{4}\zeta \left( 4 \right)\log \left( 2 \right) - \frac{9}{4}\zeta \left( 2 \right)\zeta \left( 3 \right)
\end{align}
\begin{align}\label{4.7}
\int\limits_0^1 {\frac{{\log \left( {1 + t} \right){{\log }^3}\left( t \right)}}{{1 - t}}} dt = 12\zeta \left( 5 \right) - \frac{{45}}{4}\zeta \left( 4 \right)\log \left( 2 \right) - \frac{9}{4}\zeta \left( 2 \right)\zeta \left( 3 \right),
\end{align}
\begin{align}\label{4.8}
 \int\limits_0^1 {\frac{{{{\log }^2}\left( {1 + x} \right){{\log }^2}\left( x \right)}}{{1 + x}}dx}  =& 8{\rm{L}}{{\rm{i}}_5}\left( {\frac{1}{2}} \right) + 8{\rm{L}}{{\rm{i}}_4}\left( {\frac{1}{2}} \right)\log \left( 2 \right) - \frac{{33}}{8}\zeta \left( 5 \right) - 2\zeta \left( 2 \right)\zeta \left( 3 \right)\nonumber \\
  &+ \frac{7}{2}\zeta \left( 3 \right){\log ^2}\left( 2 \right) - \frac{4}{3}\zeta \left( 2 \right){\log ^3}\left( 2 \right) + \frac{4}{{15}}{\log ^5}\left( 2 \right),
\end{align}
\begin{align}\label{4.9}
 \int\limits_0^1 {\frac{{{{\log }^2}\left( {1 + x} \right){{\log }^2}\left( {1 - x} \right)}}{x}dx}  = & 4{\rm{L}}{{\rm{i}}_5}\left( {\frac{1}{2}} \right) + 4{\rm{L}}{{\rm{i}}_4}\left( {\frac{1}{2}} \right)\log \left( 2 \right) - \frac{{25}}{8}\zeta \left( 5 \right) + \frac{7}{4}\zeta \left( 3 \right){\log ^2}\left( 2 \right)\nonumber \\
  &- \frac{2}{3}\zeta \left( 2 \right){\log ^3}\left( 2 \right) + \frac{2}{{15}}{\log ^5}\left( 2 \right),
\end{align}
\begin{align}\label{4.10}
 \int\limits_{1/2}^1 {\frac{{\log \left( x \right){{\log }^3}\left( {1 - x} \right)}}{x}dx}  = & - 6{\rm{L}}{{\rm{i}}_5}\left( {\frac{1}{2}} \right) + \frac{3}{4}\zeta \left( 4 \right)\log \left( 2 \right) - \frac{3}{8}\zeta \left( 3 \right){\log ^2}\left( 2 \right) + \frac{1}{4}{\log ^5}\left( 2 \right)\nonumber \\
  &+ 6\sum\limits_{n = 1}^\infty  {\frac{{{H_n}}}{{{n^4}{2^n}}}} ,
\end{align}
\begin{align}\label{4.11}
 \int\limits_{1/2}^1 {\frac{{{{\log }^2}\left( x \right){{\log }^2}\left( {1 - x} \right)}}{x}dx}  = & 4{\rm{L}}{{\rm{i}}_5}\left( {\frac{1}{2}} \right) + 8\zeta \left( 5 \right) - 4\zeta \left( 2 \right)\zeta \left( 3 \right) - \frac{1}{2}\zeta \left( 4 \right)\log \left( 2 \right) \nonumber \\
  &+ \frac{1}{4}\zeta \left( 3 \right){\log ^2}\left( 2 \right)+ \frac{1}{6}{\log ^5}\left( 2 \right) - 4\sum\limits_{n = 1}^\infty  {\frac{{{H_n}}}{{{n^4}{2^n}}}} ,
\end{align}
\begin{align}\label{4.12}
\int\limits_{1/2}^1 {\frac{{{{\log }^3}\left( x \right)\log \left( {1 - x} \right)}}{x}dx}  = & - 6{\rm{L}}{{\rm{i}}_5}\left( {\frac{1}{2}} \right) - 6{\rm{L}}{{\rm{i}}_4}\left( {\frac{1}{2}} \right)\log \left( 2 \right) + 6\zeta \left( 5 \right) - \frac{{21}}{8}\zeta \left( 3 \right){\log ^2}\left( 2 \right)\nonumber\\&+ \zeta \left( 2 \right){\log ^3}\left( 2 \right).
\end{align}
Next, we only prove the formulas (\ref{4.5}) and (\ref{4.9}). From Lemma \ref{lem 2.4} and \ref{lem 2.5}, we deduce that
\begin{align}\label{4.13}
 \int\limits_0^1 {\frac{{\log \left( {1 + x} \right){{\log }^m}\left( {1 - x} \right)}}{x}dx}  =& \sum\limits_{n = 1}^\infty  {\frac{{{{\left( { - 1} \right)}^{n - 1}}}}{n}} \int\limits_0^1 {{x^{n - 1}}{{\log }^m}\left( {1 - x} \right)dx} \nonumber \\
  =& {\left( { - 1} \right)^m}\sum\limits_{n = 1}^\infty  {\frac{{{Y_m}\left( n \right)}}{{{n^2}}}} {\left( { - 1} \right)^{n - 1}},
\end{align}
\begin{align}\label{4.14}
 \int\limits_0^1 {\frac{{{{\log }^2}\left( {1 + x} \right){{\log }^m}\left( {1 - x} \right)}}{x}dx}  =& 2\sum\limits_{n = 1}^\infty  {\left\{ {\frac{{{H_n}}}{n} - \frac{1}{{{n^2}}}} \right\}{{\left( { - 1} \right)}^n}} \int\limits_0^1 {{x^{n - 1}}{{\log }^m}\left( {1 - x} \right)dx} \nonumber \\
  = &2{\left( { - 1} \right)^{m + 1}}\sum\limits_{n = 1}^\infty  {\left\{ {\frac{{{H_n}{Y_m}\left( n \right)}}{{{n^2}}} - \frac{{{Y_m}\left( n \right)}}{{{n^3}}}} \right\}{{\left( { - 1} \right)}^{n - 1}}} .
\end{align}
Setting $m=2$ in (\ref{4.13}) and (\ref{4.14}) yield
\begin{align}\label{4.15}
\int\limits_0^1 {\frac{{\log \left( {1 + x} \right){{\log }^2}\left( {1 - x} \right)}}{x}dx}  = \sum\limits_{n = 1}^\infty  {\frac{{H_n^2 + H_n^{\left( 2 \right)}}}{{{n^2}}}} {\left( { - 1} \right)^{n - 1}},
\end{align}
\begin{align}\label{4.16}
\int\limits_0^1 {\frac{{{{\log }^2}\left( {1 + x} \right){{\log }^2}\left( {1 - x} \right)}}{x}dx}  = 2\sum\limits_{n = 1}^\infty  {\left\{ {\frac{{H_n^2 + H_n^{\left( 2 \right)}}}{{{n^3}}} - \frac{{H_n^3 + {H_n}H_n^{\left( 2 \right)}}}{{{n^2}}}} \right\}{{\left( { - 1} \right)}^{n - 1}}} .
\end{align}
From \cite{BBG1994,FS1998,Xyz2016,XC2016}, we know that
\begin{align*}
&\sum\limits_{n = 1}^\infty  {\frac{{H_n^2}}{{{n^2}}}} {\left( { - 1} \right)^{n - 1}} = \frac{{41}}{{16}}\zeta \left( 4 \right) + \frac{1}{2}\zeta \left( 2 \right){\log ^2}\left(2\right) - \frac{1}{{12}}{\log ^4}\left(2\right) - \frac{7}{4}\zeta \left( 3 \right)\log \left(2\right) - 2{\rm{L}}{{\rm{i}}_4}\left( {\frac{1}{2}} \right),\\
&\sum\limits_{n = 1}^\infty  {\frac{{{H^{(2)} _n}}}{{{n^2}}}{{\left( { - 1} \right)}^{n - 1}}}  =  - \frac{{51}}{{16}}\zeta (4) + 4{\rm{L}}{{\rm{i}}_4}\left( {\frac{1}{2}} \right) + \frac{7}{2}\zeta (3)\log \left(2\right) - \zeta (2){\log ^2}\left(2\right) + \frac{1}{6}{\log ^4}\left(2\right)
\end{align*}
and
\begin{align*}
\sum\limits_{n = 1}^\infty  {\frac{{H_n^{\left( 2 \right)}}}{{{n^3}}}} {\left( { - 1} \right)^{n - 1}} = & \frac{5}{8}\zeta \left( 2 \right)\zeta \left( 3 \right) - \frac{{11}}{{32}}\zeta \left( 5 \right),\\
 \sum\limits_{n = 1}^\infty  {\frac{{H_n^2}}{{{n^3}}}} {\left( { - 1} \right)^{n - 1}} =& 4{\rm{L}}{{\rm{i}}_5}\left( {\frac{1}{2}} \right) + 4\log \left(2\right) {\rm{L}}{{\rm{i}}_4}\left( {\frac{1}{2}} \right) + \frac{2}{{15}}{\log ^5}\left(2\right)  + \frac{7}{4}\zeta \left( 3 \right){\log ^2}\left(2\right) \\
  &- \frac{{19}}{{32}}\zeta \left( 5 \right) - \frac{2}{3}\zeta \left( 2 \right){\log ^3}\left(2\right) - \frac{{11}}{8}\zeta \left( 2 \right)\zeta \left( 3 \right),\\
\sum\limits_{n = 1}^\infty  {\frac{{H_n^3}}{{{n^2}}}} {\left( { - 1} \right)^{n - 1}} =& 6{\rm{L}}{{\rm{i}}_5}\left( {\frac{1}{2}} \right) + 6\log \left(2\right) {\rm{L}}{{\rm{i}}_4}\left( {\frac{1}{2}} \right) + \frac{1}{5}{\log ^5}\left(2\right)  + \frac{{21}}{8}\zeta \left( 3 \right){\log ^2}\left(2\right)  \\
 &- \frac{9}{4}\zeta \left( 5 \right) - \zeta \left( 2 \right){\log ^3}\left(2\right)  - \frac{{27}}{{16}}\zeta \left( 2 \right)\zeta \left( 3 \right),\\
\sum\limits_{n = 1}^\infty  {\frac{{{H_n}{H^{(2)}_n}}}{{{n^2}}}{{\left( { - 1} \right)}^{n - 1}}}
 =& - 4{\rm Li}{_5}\left( {\frac{1}{2}} \right) - 4\log \left(2\right) {\rm Li}{_4}\left( {\frac{1}{2}} \right) - \frac{2}{{15}}{\log ^5}\left(2\right)  - \frac{7}{4}\zeta \left( 3 \right){\log ^2}\left(2\right) \\
           & + \frac{{23}}{8}\zeta \left( 5 \right) + \frac{2}{3}\zeta \left( 2 \right){\log ^3}\left(2\right)  + \frac{{15}}{{16}}\zeta \left( 2 \right)\zeta \left( 3 \right).
\end{align*}
Hence, substituting the above identities into equations (\ref{4.15}) and (\ref{4.16}), by a direct calculation, we can obtain the desired results.
\section{Some evaluation of Euler type sums ${S_{{p_1}{p_2} \cdots {p_m},p}}\left( 1/2\right)$}
We have used our equation system to obtain explicit evaluation for all sums ${S_{{p_1}{p_2} \cdots {p_m},p}}\left( 1/2\right)$ with weight less than or equal to five. In this section, we only prove the results of all sums weight $= 5$. The formulas of weight $\leq 4$ are easily obtained.

\subsection{Weight $\leq 4$}

\begin{align*}
 &S_{1,1}\left(\frac 1{2}\right) = \frac{1}{2}\zeta \left( 2 \right), \\
 &S_{2,1}\left(\frac 1{2}\right) = \frac{5}{8}\zeta \left( 3 \right), \\
 &S_{1^2,1}\left(\frac 1{2}\right) = \frac{7}{8}\zeta \left( 3 \right), \\
&S_{1,2}\left(\frac 1{2}\right) = \zeta \left( 3 \right) - \frac{1}{2}\zeta \left( 2 \right)\log\left( 2 \right), \\
 &S_{1,3}\left(\frac 1{2}\right)  = {\rm{L}}{{\rm{i}}_4}\left( {\frac{1}{2}} \right){\rm{ + }}\frac{1}{8}\zeta \left( 4 \right) - \frac{1}{8}\zeta \left( 3 \right)\log\left( 2 \right) + \frac{1}{{24}}{\log ^4}\left( 2 \right), \\
& S_{2,2}\left(\frac 1{2}\right)= {\rm{L}}{{\rm{i}}_4}\left( {\frac{1}{2}} \right) + \frac{1}{{16}}\zeta (4) + \frac{1}{4}\zeta (3)\log\left( 2 \right) - \frac{1}{4}\zeta \left( 2 \right){\log ^2}\left( 2 \right) + \frac{1}{{24}}{\log ^4}\left( 2 \right), \\
 &S_{3,1}\left(\frac 1{2}\right) = {\rm{L}}{{\rm{i}}_4}\left( {\frac{1}{2}} \right) - \frac{5}{{16}}\zeta \left( 4 \right) + \frac{7}{8}\zeta \left( 3 \right)\log \left( 2 \right) - \frac{1}{4}\zeta \left( 2 \right){\log ^2}2\left( 2 \right) + \frac{1}{{24}}{\log ^4}\left( 2 \right), \\
 &S_{1^2,2}\left(\frac 1{2}\right) =  - {\rm{L}}{{\rm{i}}_4}\left( {\frac{1}{2}} \right)\; + \frac{{37}}{{16}}\zeta \left( 4 \right) - \frac{7}{4}\zeta \left( 3 \right)\log\left( 2 \right) + \frac{1}{4}\zeta \left( 2 \right){\log ^2}2\left( 2 \right) - \frac{1}{{24}}{\log ^4}\left( 2 \right), \\
 &S_{12,1}\left(\frac 1{2}\right)  = {\rm{L}}{{\rm{i}}_4}\left( {\frac{1}{2}} \right) - \frac{1}{8}\zeta \left( 4 \right) + \frac{7}{8}\zeta \left( 3 \right)\log\left( 2 \right) - \frac{1}{4}\zeta \left( 2 \right){\log ^2}\left( 2 \right) + \frac{1}{{24}}{\log ^4}\left( 2 \right), \\
 &S_{1^3,1}\left(\frac 1{2}\right)  =  - 5{\rm{L}}{{\rm{i}}_4}\left( {\frac{1}{2}} \right) + \frac{{25}}{4}\zeta \left( 4 \right) - \frac{{35}}{8}\zeta \left( 3 \right)\log \left( 2 \right){\rm{ + }}\frac{5}{4}\zeta \left( 2 \right){\log ^2}\left( 2 \right) - \frac{5}{{24}}{\log ^5}\left( 2 \right).
\end{align*}

\subsection{Weight $= 5$}
\begin{align}\label{5.1}
S_{1,4}\left(\frac 1{2}\right)=& 2{\rm{L}}{{\rm{i}}_5}\left( {\frac{1}{2}} \right) + {\rm{L}}{{\rm{i}}_4}\left( {\frac{1}{2}} \right)\log\left( 2 \right) + \frac{1}{{32}}\zeta \left( 5 \right)- \frac{1}{2}\zeta \left( 2 \right)\zeta \left( 3 \right)- \frac{1}{8}\zeta \left( 4 \right)\log\left( 2 \right)\nonumber \\
&   + \frac{1}{2}\zeta \left( 3 \right){\log ^2}2 - \frac{1}{6}\zeta \left( 2 \right){\log ^3}2 + \frac{1}{{40}}{\log ^5}2 ,
\end{align}
\begin{align}\label{5.2}
S_{2,3}\left(\frac 1{2}\right)=& - 2{\rm{L}}{{\rm{i}}_5}\left( {\frac{1}{2}} \right) - 3{\rm{L}}{{\rm{i}}_4}\left( {\frac{1}{2}} \right)\log \left( 2 \right) + \frac{{23}}{{64}}\zeta \left( 5 \right) + \frac{{23}}{{16}}\zeta \left( 2 \right)\zeta \left( 3 \right) - \frac{1}{{16}}\zeta \left( 4 \right)\log \left( 2 \right)\nonumber\\
& - \frac{{23}}{{16}}\zeta \left( 3 \right){\log ^2}\left( 2 \right) + \frac{7}{{12}}\zeta \left( 2 \right){\log ^3}\left( 2 \right) - \frac{{13}}{{120}}{\log ^5}\left( 2 \right),
\end{align}
\begin{align}\label{5.3}
S_{3,2}\left(\frac 1{2}\right)=&4{\rm{L}}{{\rm{i}}_5}\left( {\frac{1}{2}} \right) + 3{\rm{L}}{{\rm{i}}_4}\left( {\frac{1}{2}} \right)\log \left( 2 \right) - \frac{{81}}{{64}}\zeta \left( 5 \right) - \frac{7}{8}\zeta \left( 2 \right)\zeta \left( 3 \right) + \frac{5}{{16}}\zeta \left( 4 \right)\log \left( 2 \right)\nonumber\\
& + \frac{7}{8}\zeta \left( 3 \right){\log ^2}\left( 2 \right) - \frac{5}{{12}}\zeta \left( 2 \right){\log ^3}\left( 2 \right) + \frac{{11}}{{120}}{\log ^5}\left( 2 \right),
\end{align}
\begin{align}\label{5.4}
S_{4,1}\left(\frac 1{2}\right)=&- {\rm{L}}{{\rm{i}}_5}\left( {\frac{1}{2}} \right) - {\rm{L}}{{\rm{i}}_4}\left( {\frac{1}{2}} \right)\log \left( 2 \right) + \frac{{27}}{{32}}\zeta \left( 5 \right) + \frac{7}{{16}}\zeta \left( 2 \right)\zeta \left( 3 \right) - \frac{7}{{16}}\zeta \left( 3 \right){\log ^2}\left( 2 \right)\nonumber\\
& + \frac{1}{6}\zeta \left( 2 \right){\log ^3}\left( 2 \right) - \frac{1}{{30}}{\log ^5}\left( 2 \right),
\end{align}
\begin{align}\label{5.5}
 S_{1^2,3}\left(\frac 1{2}\right)=&- 2{\rm{L}}{{\rm{i}}_5}\left( {\frac{1}{2}} \right) - {\rm{L}}{{\rm{i}}_4}\left( {\frac{1}{2}} \right)\log \left( 2 \right) + \frac{{279}}{{64}}\zeta \left( 5 \right) - \frac{9}{{16}}\zeta \left( 2 \right)\zeta \left( 3 \right) - \frac{{37}}{{16}}\zeta \left( 4 \right)\log \left( 2 \right)\nonumber \\
  &+ \frac{7}{{16}}\zeta \left( 3 \right){\log ^2}\left( 2 \right) + \frac{1}{{12}}\zeta \left( 2 \right){\log ^3}\left( 2 \right) - \frac{1}{{40}}{\log ^5}\left( 2 \right),
 \end{align}
\begin{align}\label{5.6}
S_{1^3,2}\left(\frac 1{2}\right)=&- 14{\rm{L}}{{\rm{i}}_5}\left( {\frac{1}{2}} \right) - 9{\rm{L}}{{\rm{i}}_4}\left( {\frac{1}{2}} \right)\log \left( 2 \right) + \frac{{279}}{{16}}\zeta \left( 5 \right) - \frac{7}{8}\zeta \left( 2 \right)\zeta \left( 3 \right) - \frac{{25}}{4}\zeta \left( 4 \right)\log \left( 2 \right)\nonumber\\
& - \frac{7}{4}\zeta \left( 3 \right){\log ^2}\left( 2 \right) + \frac{{13}}{{12}}\zeta \left( 2 \right){\log ^3}\left( 2 \right) - \frac{{31}}{{120}}{\log ^5}\left( 2 \right),
\end{align}
\begin{align}\label{5.7}
S_{12,2}\left(\frac 1{2}\right)=&2{\rm{L}}{{\rm{i}}_5}\left( {\frac{1}{2}} \right) + {\rm{L}}{{\rm{i}}_4}\left( {\frac{1}{2}} \right)\log \left( 2 \right) - \frac{{31}}{{32}}\zeta \left( 5 \right) + \frac{1}{8}\zeta \left( 2 \right)\zeta \left( 3 \right) + \frac{1}{8}\zeta \left( 4 \right)\log \left( 2 \right)\nonumber\\
& - \frac{1}{{12}}\zeta \left( 2 \right){\log ^3}\left( 2 \right) + \frac{1}{{40}}{\log ^5}\left( 2 \right),
\end{align}
\begin{align}\label{5.8}
S_{13,1}\left(\frac 1{2}\right)=&3{\rm{L}}{{\rm{i}}_5}\left( {\frac{1}{2}} \right) + 3{\rm{L}}{{\rm{i}}_4}\left( {\frac{1}{2}} \right)\log \left( 2 \right) - \frac{{31}}{{64}}\zeta \left( 5 \right) - \frac{7}{8}\zeta \left( 2 \right)\zeta \left( 3 \right) + \frac{{21}}{{16}}\zeta \left( 3 \right){\log ^2}\left( 2 \right)\nonumber\\
& - \frac{1}{2}\zeta \left( 2 \right){\log ^3}\left( 2 \right) + \frac{1}{{10}}{\log ^5}\left( 2 \right),
\end{align}
\begin{align}\label{5.9}
S_{1^22,1}\left(\frac 1{2}\right)=& 3{\rm{L}}{{\rm{i}}_5}\left( {\frac{1}{2}} \right) + 3{\rm{L}}{{\rm{i}}_4}\left( {\frac{1}{2}} \right)\log \left( 2 \right) - \frac{{31}}{{32}}\zeta \left( 5 \right) - \frac{7}{{16}}\zeta \left( 2 \right)\zeta \left( 3 \right) + \frac{{21}}{{16}}\zeta \left( 3 \right){\log ^2}\left( 2 \right)\nonumber\\
& - \frac{1}{2}\zeta \left( 2 \right){\log ^3}\left( 2 \right) + \frac{1}{{10}}{\log ^5}\left( 2 \right),
\end{align}
\begin{align}\label{5.10}
S_{1^4,1}\left(\frac 1{2}\right)=&- 15{\rm{L}}{{\rm{i}}_5}\left( {\frac{1}{2}} \right) - 15{\rm{L}}{{\rm{i}}_4}\left( {\frac{1}{2}} \right)\log \left( 2 \right) + \frac{{341}}{{16}}\zeta \left( 5 \right) - \frac{{35}}{{16}}\zeta \left( 2 \right)\zeta \left( 3 \right) - \frac{{105}}{{16}}\zeta \left( 3 \right){\log ^2}\left( 2 \right)\nonumber\\
& + \frac{5}{2}\zeta \left( 2 \right){\log ^3}\left( 2 \right) - \frac{1}{2}{\log ^5}\left( 2 \right),
\end{align}
\begin{align}\label{5.11}
S_{2^2,1}\left(\frac 1{2}\right)=&- 7{\rm{L}}{{\rm{i}}_5}\left( {\frac{1}{2}} \right) - 7{\rm{L}}{{\rm{i}}_4}\left( {\frac{1}{2}} \right)\log \left( 2 \right) + \frac{{31}}{{16}}\zeta \left( 5 \right){\rm{ + }}\frac{{49}}{{16}}\zeta \left( 2 \right)\zeta \left( 3 \right) - \frac{{49}}{{16}}\zeta \left( 3 \right){\log ^2}\left( 2 \right)\nonumber\\
& + \frac{7}{6}\zeta \left( 2 \right){\log ^3}\left( 2 \right) - \frac{7}{{30}}{\log ^5}\left( 2 \right).
\end{align}

\subsection{Weight $= 6$}
\begin{align}\label{5.12}
S_{1,5}\left(\frac1{2}\right) =& 3{\rm{L}}{{\rm{i}}_6}\left( {\frac{1}{2}} \right) + {\rm{L}}{{\rm{i}}_5}\left( {\frac{1}{2}} \right)\log \left( 2 \right) - \frac{1}{2}\zeta \left( {\bar 5,1} \right) - \frac{{51}}{{32}}\zeta \left( 6 \right) - \frac{1}{4}{\zeta ^2}\left( 3 \right) - \frac{1}{{32}}\zeta \left( 5 \right)\log \left( 2 \right)\nonumber\\
&+ \frac{1}{2}\zeta \left( 2 \right)\zeta \left( 3 \right)\log \left( 2 \right) + \frac{1}{{16}}\zeta \left( 4 \right){\log ^2}\left( 2 \right) - \frac{1}{6}\zeta \left( 3 \right){\log ^3}\left( 2 \right)\nonumber\\& + \frac{1}{{24}}\zeta \left( 2 \right){\log ^4}\left( 2 \right) - \frac{1}{{240}}{\log ^6}\left( 2 \right).
\end{align}

\subsection{Proof of all sums of weight $=5$}
In \cite{FS1998}, Flajolet and Salvy gave an explicit formula for alternating Euler sums ${\bar S}_{1,m}:=S_{1,m}\left(-1\right)$ in term of zeta values, polylogarithms and $\log\left(2\right)$ when $m$ is a even by using the method of contour integral representations and residue computation. Hence, we deduce the result
\begin{align}\label{5.13}
\sum\limits_{n{\rm{ = }}1}^\infty  {\frac{{{H_n}}}{{{n^4}}}{{\left( { - 1} \right)}^{n - 1}} = } \frac{{59}}{{32}}\zeta \left( 5 \right) - \frac{1}{2}\zeta \left( 2 \right)\zeta \left( 3 \right).
\end{align}
Letting $m=3$ in (\ref{3.1}) and combining formulas (\ref{4.10})-(\ref{4.12}), we have
\begin{align}\label{5.14}
\sum\limits_{n = 1}^\infty  {\frac{{{H_n}}}{{{n^4}}}} {\left( { - 1} \right)^n} = &3{S_{1,4}}\left( {\frac{1}{2}} \right) - 6{\rm{L}}{{\rm{i}}_5}\left( {\frac{1}{2}} \right) - 3{\rm{L}}{{\rm{i}}_4}\left( {\frac{1}{2}} \right)\log \left( 2 \right) - \frac{{31}}{{16}}\zeta \left( 5 \right) + 2\zeta \left( 2 \right)\zeta \left( 3 \right)\nonumber\\
&+ \frac{3}{8}\zeta \left( 4 \right)\log \left( 2 \right) - \frac{3}{2}\zeta \left( 3 \right){\log ^2}\left( 2 \right) + \frac{1}{2}\zeta \left( 2 \right){\log ^3}\left( 2 \right) - \frac{3}{{40}}{\log ^5}\left( 2 \right).
\end{align}
Thus, applying (\ref{5.13}) to (\ref{5.14}), the result is (\ref{5.1}).

Next, we prove the identities (\ref{5.2})-(\ref{5.4}). Multiplying (\ref{3.9}) by $\frac {\log\left(1-x\right)}{x}$ and integrating over the interval $(0,1)$, and using (\ref{2.7}), we obtain
\begin{align}\label{5.15}
&\sum\limits_{n = 1}^\infty  {\frac{{{H_n}L_n^2\left( 1 \right)}}{{{n^2}}}}  - 2\sum\limits_{n = 1}^\infty  {\frac{{{H_n}{L_n}\left( 1 \right)}}{{{n^3}}}{{\left( { - 1} \right)}^{n - 1}}}  + 2\sum\limits_{n = 1}^\infty  {\frac{{{H_n}}}{{{n^4}}}}  - \sum\limits_{n = 1}^\infty  {\frac{{{H_n}H_n^{\left( 2 \right)}}}{{{n^2}}}} \nonumber\\
& = \int\limits_0^1 {\frac{{4\left( {{\rm{L}}{{\rm{i}}_3}\left( {\frac{{1 - x}}{2}} \right) - {\rm{L}}{{\rm{i}}_3}\left( {\frac{1}{2}} \right)} \right) - 2\log \left( {1 - x} \right)\left( {{\rm{L}}{{\rm{i}}_2}\left( {\frac{{1 - x}}{2}} \right) + {\rm{L}}{{\rm{i}}_2}\left( {\frac{1}{2}} \right)} \right)}}{x}} \log \left( {1 - x} \right)dx \nonumber\\
&\mathop  = \limits^{x = 1 - t} \int\limits_0^1 {\frac{{4\left( {{\rm{L}}{{\rm{i}}_3}\left( {\frac{t}{2}} \right) - {\rm{L}}{{\rm{i}}_3}\left( {\frac{1}{2}} \right)} \right) - 2\log \left( t \right)\left( {{\rm{L}}{{\rm{i}}_2}\left( {\frac{t}{2}} \right) + {\rm{L}}{{\rm{i}}_2}\left( {\frac{1}{2}} \right)} \right)}}{{1 - t}}} \log \left( t \right)dt \nonumber\\
& = 4\int\limits_0^1 {\frac{{\log \left( t \right){\rm{L}}{{\rm{i}}_3}\left( {\frac{t}{2}} \right)}}{{1 - t}}} dt - 4{\rm{L}}{{\rm{i}}_3}\left( {\frac{1}{2}} \right)\int\limits_0^1 {\frac{{\log \left( t \right)}}{{1 - t}}} dt \nonumber\\&\quad - 2{\rm{L}}{{\rm{i}}_2}\left( {\frac{1}{2}} \right)\int\limits_0^1 {\frac{{{{\log }^2}\left( t \right)}}{{1 - t}}} dt - 2\int\limits_0^1 {\frac{{{{\log }^2}\left( t \right){\rm{L}}{{\rm{i}}_2}\left( {\frac{t}{2}} \right)}}{{1 - t}}} dt.
\end{align}
Furthermore, by using (\ref{3.17}) and (\ref{4.1}), then the integrals on the right hand side of (\ref{5.15}) can be rewritten as
\begin{align}\label{5.16}
& \int\limits_0^1 {\frac{{\log \left( t \right)}}{{1 - t}}} dt =  - \zeta \left( 2 \right),
\end{align}
\begin{align}\label{5.17}
 &\int\limits_0^1 {\frac{{{{\log }^2}\left( t \right)}}{{1 - t}}} dt = 2\zeta \left( 3 \right),
 \end{align}
\begin{align}\label{5.18}
 &\int\limits_0^1 {\frac{{\log \left( t \right){\rm{L}}{{\rm{i}}_3}\left( {\frac{t}{2}} \right)}}{{1 - t}}} dt = {S_{2,3}}\left( {\frac{1}{2}} \right) - \zeta \left( 2 \right){\rm{L}}{{\rm{i}}_3}\left( {\frac{1}{2}} \right),
\end{align}
\begin{align}\label{5.19}
 &\int\limits_0^1 {\frac{{{{\log }^2}\left( t \right){\rm{L}}{{\rm{i}}_2}\left( {\frac{t}{2}} \right)}}{{1 - t}}}  = 2\zeta \left( 3 \right){\rm{L}}{{\rm{i}}_2}\left( {\frac{1}{2}} \right) - 2{S_{3,2}}\left( {\frac{1}{2}} \right).
\end{align}
Therefore, combining formulas (\ref{5.15})-(\ref{5.19}), by a simple calculation, we obtain the following equation
\begin{align}\label{5.20}
\sum\limits_{n = 1}^\infty  {\frac{{{H_n}L_n^2\left( 1 \right)}}{{{n^2}}}}  - 2\sum\limits_{n = 1}^\infty  {\frac{{{H_n}{L_n}\left( 1 \right)}}{{{n^3}}}{{\left( { - 1} \right)}^{n - 1}}}  =& 4{S_{2,3}}\left( {\frac{1}{2}} \right) + 4{S_{3,2}}\left( {\frac{1}{2}} \right) + 4\zeta \left( 3 \right){\log ^2}\left( 2 \right)\nonumber\\& - 5\zeta \left( 5 \right) - \zeta \left( 2 \right)\zeta \left( 3 \right).
\end{align}
On the other hand, from (3.13) of \cite{Xu2015}, we deduce that
\begin{align}\label{5.21}
&\sum\limits_{n = 1}^\infty  {\frac{{{H_n}L_n^2\left( 1 \right) + {H_n}H_n^{\left( 2 \right)} + H_n^2{L_n}\left( 1 \right) + {L_n}\left( 1 \right){L_n}\left( 2 \right)}}{{{n^2}}}}\nonumber \\
&= 2\zeta \left( 2 \right)\zeta \left( 3 \right) + 2\log \left( 2 \right)\sum\limits_{n = 1}^\infty  {\frac{{H_n^2 + {H_n}{L_n}\left( 1 \right)}}{{{n^2}}}}  - 2\log \left( 2 \right)\sum\limits_{n = 1}^\infty  {\frac{{{H_n}}}{{{n^3}}}\left( {1 + {{\left( { - 1} \right)}^{n - 1}}} \right)}\nonumber \\
&\quad + 2\sum\limits_{n = 1}^\infty  {\frac{{{H_n}{L_n}\left( 1 \right)}}{{{n^3}}}{{\left( { - 1} \right)}^{n - 1}}}  - 2\zeta \left( 2 \right)\sum\limits_{n = 1}^\infty  {\frac{{{L_n}\left( 1 \right)}}{{{n^2}}}}  + 2\sum\limits_{n = 1}^\infty  {\frac{{{H_n}{L_n}\left( 1 \right)}}{{{n^3}}}} .
\end{align}
Hence, substituting (\ref{5.20}) into (\ref{5.21}) with the help of the results of alternating Euler sums in the reference \cite{XC2016}, we can get the following relation
\begin{align}\label{5.22}
 {S_{2,3}}\left( {\frac{1}{2}} \right) + {S_{3,2}}\left( {\frac{1}{2}} \right) = &2{\rm{L}}{{\rm{i}}_5}\left( {\frac{1}{2}} \right) - \frac{{29}}{{32}}\zeta \left( 5 \right) + \frac{9}{{16}}\zeta \left( 2 \right)\zeta \left( 3 \right) + \frac{1}{4}\zeta \left( 4 \right)\log \left( 2 \right)\nonumber \\
  &- \frac{9}{{16}}\zeta \left( 3 \right){\log ^2}\left( 2 \right) + \frac{1}{6}\zeta \left( 2 \right){\log ^3}\left( 2 \right) - \frac{1}{{60}}{\log ^5}\left( 2 \right).
\end{align}
Taking $(s,t)=(2,3)$ and $(1,4)$ in (\ref{2.1}), then letting $x=\frac1{2}$, we have
\begin{align}\label{5.23}
3{S_{2,3}}\left( {\frac{1}{2}} \right) + {S_{3,2}}\left( {\frac{1}{2}} \right) =&  - 2{\rm{L}}{{\rm{i}}_5}\left( {\frac{1}{2}} \right) - 6{\rm{L}}{{\rm{i}}_4}\left( {\frac{1}{2}} \right)\log \left( 2 \right) - \frac{3}{{16}}\zeta \left( 5 \right) + \frac{{55}}{{16}}\zeta \left( 2 \right)\zeta \left( 3 \right)+ \frac{1}{8}\zeta \left( 4 \right)\log \left( 2 \right) \nonumber\\
&  - \frac{{55}}{{16}}\zeta \left( 3 \right){\log ^2}\left( 2 \right) + \frac{4}{3}\zeta \left( 2 \right){\log ^3}\left( 2 \right) - \frac{7}{{30}}{\log ^5}\left( 2 \right),
\end{align}
\begin{align}\label{5.24}
2{S_{1,4}}\left( {\frac{1}{2}} \right) + {S_{2,3}}\left( {\frac{1}{2}} \right) + {S_{3,2}}\left( {\frac{1}{2}} \right) + {S_{4,1}}\left( {\frac{1}{2}} \right) = 5{\rm{L}}{{\rm{i}}_5}\left( {\frac{1}{2}} \right) + {\rm{L}}{{\rm{i}}_4}\left( {\frac{1}{2}} \right)\log \left( 2 \right).
\end{align}
Combining equations (\ref{5.22})-(\ref{5.24}) we can prove the formulas (\ref{5.2})-(\ref{5.4}).

Now, we establish enough equations of Euler type sums to prove the identities (\ref{5.5})-(\ref{5.11}). In (\ref{2.2}), letting $m=k=2$, then
\begin{align}\label{5.25}
&2{\log ^2}\left( 2 \right)\zeta \left( {2,1;\frac{1}{2}} \right) + 4\log \left( 2 \right)\zeta \left( {3,1;\frac{1}{2}} \right) + 4\zeta \left( {4,1;\frac{1}{2}} \right)\nonumber\\
& + 2{\log ^2}\left( 2 \right)\zeta \left( {{{\left\{ 1 \right\}}_3};\frac{1}{2}} \right) + 4\log \left( 2 \right)\zeta \left( {2,1,1;\frac{1}{2}} \right) + 4\zeta \left( {3,1,1;\frac{1}{2}} \right)\nonumber\\
& = 4\zeta \left( {3,1,1} \right).
\end{align}
By the definition of Euler type sums ${S_{{p_1}{p_2} \cdots {p_m},p}}\left( x \right)$ and multiple polylogarithm function $\zeta \left( {{s_1},{s_2}, \cdots ,{s_m};x} \right)$, it is easily seen that
\[{S_{p,q}}\left( {\frac{1}{2}} \right) = \zeta \left( {q,p;\frac{1}{2}} \right) + {\rm{L}}{{\rm{i}}_{p + q}}\left( {\frac{1}{2}} \right).\]
Noting that the multiple zeta value $\zeta(m+1,\{1\}_{k-1} )$ can be
represented as a polynomial of zeta values with rational coefficients (see \cite{BBBL1997,Xu2017}),
and using (\ref{2.2}), we arrive at the conclusion that
\begin{align}\label{5.26}
 \zeta \left( {3,1,1;\frac{1}{2}} \right) = & - {\rm{L}}{{\rm{i}}_5}\left( {\frac{1}{2}} \right) + \frac{{63}}{{32}}\zeta \left( 5 \right) - \frac{1}{2}\zeta \left( 2 \right)\zeta \left( 3 \right) - \zeta \left( 4 \right)\log \left( 2 \right)\nonumber \\
  &+ \frac{7}{{16}}\zeta \left( 3 \right){\log ^2}\left( 2 \right) - \frac{1}{{12}}\zeta \left( 2 \right){\log ^3}\left( 2 \right) + \frac{1}{{60}}{\log ^5}\left( 2 \right).
\end{align}
Similarly, putting $m=3$ in Lemma \ref{lem 2.3}, we conclude that
\begin{align}\label{5.27}
\zeta \left( {2,{{\left\{ 1 \right\}}_3};\frac{1}{2}} \right) =&  - {\rm{L}}{{\rm{i}}_5}\left( {\frac{1}{2}} \right) - \log \left( 2 \right){\rm{L}}{{\rm{i}}_4}\left( {\frac{1}{2}} \right) + \zeta \left( 5 \right) - \frac{7}{{16}}\zeta \left( 3 \right){\log ^2}\left( 2 \right)\nonumber\\& + \frac{1}{6}\zeta \left( 2 \right){\log ^3}\left( 2 \right) - \frac{1}{{24}}{\log ^5}\left( 2 \right).
\end{align}
Furthermore, from the definition of $\zeta \left( {{s_1},{s_2}, \cdots ,{s_m};x} \right)$, we know that
\begin{align}\label{5.28}
\zeta \left( {3,1,1;\frac{1}{2}} \right) =& \sum\limits_{n = 1}^\infty  {\frac{1}{{{n^3}{2^n}}}{\zeta _{n - 1}}\left( {1,1} \right)}\nonumber \\
 =& \sum\limits_{n = 1}^\infty  {\frac{1}{{{n^3}{2^n}}}\left\{ {\frac{{H_n^2 - H_n^{\left( 2 \right)}}}{2} - \frac{{{H_n}}}{n} + \frac{1}{{{n^2}}}} \right\}}\nonumber \\
 = &\frac{1}{2}{S_{{1^2},3}}\left( {\frac{1}{2}} \right) - \frac{1}{2}{S_{2,3}}\left( {\frac{1}{2}} \right) - {S_{1,4}}\left( {\frac{1}{2}} \right) + {\rm{L}}{{\rm{i}}_5}\left( {\frac{1}{2}} \right),
\end{align}
\begin{align}\label{5.29}
\zeta \left( {2,{{\left\{ 1 \right\}}_3};\frac{1}{2}} \right) &= \sum\limits_{n = 1}^\infty  {\frac{1}{{{n^2}{2^n}}}{\zeta _{n - 1}}\left( {{{\left\{ 1 \right\}}_3}} \right)}\nonumber \\
& = \frac{1}{6}\sum\limits_{n = 1}^\infty  {\frac{1}{{{n^2}{2^n}}}\left\{ {H_n^3 - 3{H_n}H_n^{\left( 2 \right)} + 2H_n^{\left( 3 \right)} - 3\frac{{H_n^2 - H_n^{\left( 2 \right)}}}{n} + 6\frac{{{H_n}}}{{{n^2}}} - 6\frac{1}{{{n^3}}}} \right\}}\nonumber \\
  &= \frac{1}{6}\left[ {{S_{{1^3},2}}\left( {\frac{1}{2}} \right) - 3{S_{12,2}}\left( {\frac{1}{2}} \right) + 2{S_{3,2}}\left( {\frac{1}{2}} \right)} \right]\nonumber \\
  &\quad- \frac{1}{2}\left[ {{S_{{1^2},3}}\left( {\frac{1}{2}} \right) - {S_{2,3}}\left( {\frac{1}{2}} \right)} \right] + {S_{1,4}}\left( {\frac{1}{2}} \right) - {\rm{L}}{{\rm{i}}_5}\left( {\frac{1}{2}} \right).
\end{align}
Setting $m=1,k=3$ in (\ref{3.14}) with the help of (\ref{4.7}) we obtain
\begin{align}\label{5.30}
{S_{{1^3},2}}\left( {\frac{1}{2}} \right) + 3{S_{12,2}}\left( {\frac{1}{2}} \right) = & - 8{\rm{L}}{{\rm{i}}_5}\left( {\frac{1}{2}} \right) - 6\log \left( 2 \right){\rm{L}}{{\rm{i}}_4}\left( {\frac{1}{2}} \right) + \frac{{465}}{{32}}\zeta \left( 5 \right) - \frac{1}{2}\zeta \left( 2 \right)\zeta \left( 3 \right)\nonumber\\
& - \frac{{47}}{8}\zeta \left( 4 \right)\log \left( 2 \right) - \frac{7}{4}\zeta \left( 3 \right){\log ^2}\left( 2 \right) + \frac{5}{6}\zeta \left( 2 \right){\log ^3}\left( 2 \right)\nonumber\\& - \frac{{11}}{{60}}{\log ^5}\left( 2 \right).
\end{align}
Thus, the relations (\ref{5.26})-(\ref{5.30}) yield the results (\ref{5.5})-(\ref{5.7}). Then, by using formulas (\ref{2.5}), (\ref{2.6}) and the recurrence relation of Stirling numbers of the first kind, we can find that
\begin{align}\label{5.31}
\sum\limits_{n = 1}^\infty  {\frac{{s\left( {n + 1,k} \right)}}{{n!{n^m}{2^n}}} = \zeta \left( {m + 1,{{\left\{ 1 \right\}}_{k - 2}};\frac{1}{2}} \right)}  + \zeta \left( {m,{{\left\{ 1 \right\}}_{k - 1}};\frac{1}{2}} \right),\;k \ge 2,m \ge 1.
\end{align}
Letting $k=5,m=1$ in above equation we get the relation
\begin{align}\label{5.32}
 &\sum\limits_{n = 1}^\infty  {\frac{{H_n^4 - 6H_n^2H_n^{\left( 2 \right)} + 8{H_n}H_n^{\left( 3 \right)} + 3{{\left( {H_n^{\left( 2 \right)}} \right)}^2} - 6H_n^{\left( 4 \right)}}}{{n{2^n}}}}\nonumber  \\
  = & - 24{\rm{L}}{{\rm{i}}_5}\left( {\frac{1}{2}} \right) - 24\log \left( 2 \right){\rm{L}}{{\rm{i}}_4}\left( {\frac{1}{2}} \right) + 24\zeta \left( 5 \right)\nonumber \\
  &- \frac{{21}}{2}\zeta \left( 3 \right){\log ^2}\left( 2 \right) + 4\zeta \left( 2 \right){\log ^3}\left( 2 \right) - \frac{4}{5}{\log ^5}\left( 2 \right).
\end{align}
In Theorem \ref{thm3.4}, taking $(m,k)=(0,4)$ and $(1,3)$ with the help of formulas (\ref{4.2}) and
\[\int\limits_0^1 {\frac{{\log \left( {1 + t} \right){{\log }^3}\left( t \right)}}{{1 + t}}dt}  = \frac{{87}}{{16}}\zeta \left( 5 \right) - 3\zeta \left( 2 \right)\zeta \left( 3 \right),\]
we obtain the following two equations
\begin{align}\label{5.33}
\sum\limits_{n = 1}^\infty  {\frac{{H_n^4 + 6H_n^2H_n^{\left( 2 \right)} + 8{H_n}H_n^{\left( 3 \right)} + 3{{\left( {H_n^{\left( 2 \right)}} \right)}^2} + 6H_n^{\left( 4 \right)}}}{{n{2^n}}}}  = \frac{{45}}{2}\zeta \left( 5 \right),
\end{align}
\begin{align}\label{5.34}
\sum\limits_{n = 1}^\infty  {\frac{{H_n^4 + 3H_n^2H_n^{\left( 2 \right)} + 2{H_n}H_n^{\left( 3 \right)}}}{{n{2^n}}}}  = \frac{{279}}{{16}}\zeta \left( 5 \right) - \frac{{21}}{4}\zeta \left( 2 \right)\zeta \left( 3 \right).
\end{align}
Futhermore, setting $m=2$ in (\ref{3.16}) and using the formulas (\ref{4.5}) and (\ref{4.9}), by a simple calculation, we deduce the identity (\ref{5.8}), namely
\begin{align}\label{5.35}
\sum\limits_{n = 1}^\infty  {\frac{{{H_n}H_n^{\left( 3 \right)}}}{{n{2^n}}}}  =& \sum\limits_{n = 1}^\infty  {\frac{{H_n^{\left( 3 \right)}}}{{{n^2}{2^n}}}}  - \frac{1}{4}\int\limits_0^1 {\frac{{{{\log }^2}\left( {1 + t} \right){{\log }^2}\left( {1 - t} \right)}}{t}dt}\nonumber \\
& + \frac{1}{2}\log \left( 2 \right)\int\limits_0^1 {\frac{{\log \left( {1 + t} \right){{\log }^2}\left( {1 - t} \right)}}{t}dt}\nonumber \\
  =& 3{\rm{L}}{{\rm{i}}_5}\left( {\frac{1}{2}} \right) + 3\log \left( 2 \right){\rm{L}}{{\rm{i}}_4}\left( {\frac{1}{2}} \right) - \frac{{31}}{{64}}\zeta \left( 5 \right) - \frac{7}{8}\zeta \left( 2 \right)\zeta \left( 3 \right) \nonumber\\
 &+ \frac{{21}}{{16}}\zeta \left( 3 \right){\log ^2}\left( 2 \right) - \frac{1}{2}\zeta \left( 2 \right){\log ^3}\left( 2 \right) + \frac{1}{{10}}{\log ^5}\left( 2 \right).
\end{align}
Finally, combining (\ref{5.32})-(\ref{5.35}) we can prove the identities (\ref{5.9}) and (\ref{5.11}). The proofs of formulas (\ref{5.1})-(\ref{5.11}) are finished.

Letting $m=4$ in (\ref{3.1}), by a similar argument as in the proof of (\ref{5.1}), we have the formula (\ref{5.12}).

Moreover, we can evaluate the following two alternating Euler sums
\begin{align*}
&\sum\limits_{n = 1}^\infty  {\frac{{{H_n}L_n^2\left( 1 \right)}}{n}} {\left( { - 1} \right)^{n - 1}} = \frac{7}{4}\zeta \left( 2 \right){\log ^2}\left( 2 \right)  - \frac{1}{4}{\log ^4}\left( 2 \right)  + \frac{3}{8}\zeta (3)\log \left( 2 \right),\\
&\sum\limits_{n = 1}^\infty  {\frac{{{H_n}L_n^2\left( 1 \right)}}{{{n^2}}}}  = 12{\rm{L}}{{\rm{i}}_5}\left( {\frac{1}{2}} \right) - \frac{{53}}{4}\zeta \left( 5 \right) + \zeta \left( 2 \right)\zeta \left( 3 \right) + 9\zeta \left( 4 \right)\log \left( 2 \right) + \zeta \left( 2 \right){\log ^3}\left( 2 \right) - \frac{1}{{10}}{\log ^5}\left( 2 \right).
\end{align*}
The proofs of above identities are left to the readers.
\section{Conclusion}

In this paper, we have proved the conclusion: all Euler type sums ${S_{{p_1}{p_2} \cdots {p_m},p}}\left( 1/2\right)$  of weight $\leq 5$  are reducible to Q-linear combinations of single zeta values, polylogarithms and $\log(2)$. Based on the above discussion, we conjectured that all such sums with $w=6$ satisfies a relation involving homogeneous combinations of these constants
\begin{align*}
&{\rm{L}}{{\rm{i}}_6}\left( {\frac{1}{2}} \right),{\rm{L}}{{\rm{i}}_5}\left( {\frac{1}{2}} \right)\log \left( 2 \right),{\rm{L}}{{\rm{i}}_4}\left( {\frac{1}{2}} \right){\log ^2}\left( 2 \right),{\rm{L}}{{\rm{i}}_4}\left( {\frac{1}{2}} \right)\zeta \left( 2 \right),\\
&\zeta \left( 6 \right),{\zeta ^2}\left( 3 \right),\zeta \left( 5 \right)\log \left( 2 \right),\zeta \left( 2 \right)\zeta \left( 3 \right)\log \left( 2 \right),\zeta \left( 4 \right){\log ^2}\left( 2 \right),\\
&\zeta \left( 3 \right){\log ^3}\left( 2 \right),\zeta \left( 2 \right){\log ^4}\left( 2 \right),{\log ^6}\left( 2 \right).
\end{align*}
However, we have been unable, so far, to prove the conjecture. By using the method of this paper, we can establish some identities involving two or more Euler type sums of the weight $=6$. Some these relations are shown in following  
\begin{align}\label{6.1}
 {\rm{Li}}_3^2\left( {\frac{1}{2}} \right) = 2{S_{3,3}}\left( {\frac{1}{2}} \right) + 6{S_{2,4}}\left( {\frac{1}{2}} \right) + 12{S_{1,5}}\left( {\frac{1}{2}} \right) - 20{\rm{L}}{{\rm{i}}_6}\left( {\frac{1}{2}} \right), 
\end{align}
\begin{align}\label{6.2}
 {\rm{L}}{{\rm{i}}_2}\left( {\frac{1}{2}} \right){\rm{L}}{{\rm{i}}_4}\left( {\frac{1}{2}} \right) = {S_{4,2}}\left( {\frac{1}{2}} \right) + 2{S_{3,3}}\left( {\frac{1}{2}} \right) + 4{S_{2,4}}\left( {\frac{1}{2}} \right) + 8{S_{1,5}}\left( {\frac{1}{2}} \right) - 15{\rm{L}}{{\rm{i}}_6}\left( {\frac{1}{2}} \right), 
\end{align}
\begin{align}\label{6.3}
 {\rm{L}}{{\rm{i}}_5}\left( {\frac{1}{2}} \right)\log \left( 2 \right) = {S_{5,1}}\left( {\frac{1}{2}} \right) + {S_{4,2}}\left( {\frac{1}{2}} \right) + {S_{3,3}}\left( {\frac{1}{2}} \right) + {S_{2,4}}\left( {\frac{1}{2}} \right) + 2{S_{1,5}}\left( {\frac{1}{2}} \right) - 6{\rm{L}}{{\rm{i}}_6}\left( {\frac{1}{2}} \right), 
\end{align}
\begin{align}\label{6.4}
 \zeta \left( {4,1,1;\frac{1}{2}} \right) = &\frac{1}{2}{S_{{1^2},4}}\left( {\frac{1}{2}} \right) - \frac{1}{2}{S_{2,4}}\left( {\frac{1}{2}} \right) - {S_{1,5}}\left( {\frac{1}{2}} \right) + {\rm{L}}{{\rm{i}}_6}\left( {\frac{1}{2}} \right)\nonumber \\
  = &\frac{{23}}{{32}}\zeta \left( 6 \right) - \frac{1}{2}{\zeta ^2}\left( 3 \right){\rm{ + L}}{{\rm{i}}_5}\left( {\frac{1}{2}} \right)\log \left( 2 \right) + \frac{1}{2}{\rm{L}}{{\rm{i}}_4}\left( {\frac{1}{2}} \right){\log ^2}\left( 2 \right) - \frac{{63}}{{32}}\zeta \left( 5 \right)\log \left( 2 \right)\nonumber \\
  &+ \frac{1}{2}\zeta \left( 2 \right)\zeta \left( 3 \right)\log \left( 2 \right) + \frac{1}{2}\zeta \left( 4 \right){\log ^2}\left( 2 \right) - \frac{1}{{24}}\zeta \left( 2 \right){\log ^4}\left( 2 \right) + \frac{1}{{90}}{\log ^6}\left( 2 \right). 
\end{align}
From formulas (\ref{6.1})-(\ref{6.3}) we obtain the closed form of Euler type sum ${S_{5,1}}\left( {\frac{1}{2}} \right)$ follows 
\begin{align}\label{6.5}
 {S_{5,1}}\left( {\frac{1}{2}} \right) =& {\rm{L}}{{\rm{i}}_6}\left( {\frac{1}{2}} \right) + {\rm{L}}{{\rm{i}}_5}\left( {\frac{1}{2}} \right)\log \left( 2 \right) - {\rm{L}}{{\rm{i}}_2}\left( {\frac{1}{2}} \right){\rm{L}}{{\rm{i}}_4}\left( {\frac{1}{2}} \right) + \frac{1}{2}{\rm{Li}}_3^2\left( {\frac{1}{2}} \right)\nonumber \\
  =& {\rm{L}}{{\rm{i}}_6}\left( {\frac{1}{2}} \right) + {\rm{L}}{{\rm{i}}_5}\left( {\frac{1}{2}} \right)\log \left( 2 \right) - \frac{1}{2}{\rm{L}}{{\rm{i}}_4}\left( {\frac{1}{2}} \right)\zeta \left( 2 \right) + \frac{1}{2}{\rm{L}}{{\rm{i}}_4}\left( {\frac{1}{2}} \right){\log ^2}\left( 2 \right)\nonumber \\
  &+ \frac{{49}}{{128}}{\zeta ^2}\left( 3 \right) - \frac{7}{{16}}\zeta \left( 2 \right)\zeta \left( 3 \right)\log \left( 2 \right) + \frac{5}{{16}}\zeta \left( 4 \right){\log ^2}\left( 2 \right)\nonumber \\
  &+ \frac{7}{{48}}\zeta \left( 3 \right){\log ^3}\left( 2 \right) - \frac{1}{{12}}\zeta \left( 2 \right){\log ^4}\left( 2 \right) + \frac{1}{{72}}{\log ^6}\left( 2 \right).
\end{align}
More general, we have (see the Corollary 1 in the reference \cite{SAZ2012})
\[{S_{2n - 1,1}}\left( z \right) = {\rm{L}}{{\rm{i}}_{2n}}\left( z \right) + \frac{1}{2}\sum\limits_{k = 1}^{2n - 1} {{{\left( { - 1} \right)}^{k + 1}}{\rm{L}}{{\rm{i}}_k}\left( z \right)} {\rm{L}}{{\rm{i}}_{2n - k}}\left( z \right),\]
Where $n\in\N, z\in D:=\mathbb{C}\setminus\{z:\left| {\arg \left( {1 - z} \right)} \right| < \pi \}$.

\end{document}